\documentclass[12pt,draftcls,onecolumn]{IEEEtran}

\usepackage{amsfonts}
\usepackage{amssymb}
\usepackage{amsmath}
\usepackage{mathptmx}
\usepackage{cancel}
\usepackage{pstricks}
\usepackage{psfrag}
\usepackage{mathrsfs}
\usepackage{graphicx}
\usepackage{pgf,tikz}
\usepackage{algorithm}
\usepackage{algpseudocode}
\usepackage{soul}
\usetikzlibrary{arrows}
\def\be{\begin{equation}}
\def\ee{\end{equation}}
\def\bea{\begin{eqnarray}}
\def\eea{\end{eqnarray}}
\def\beann{\begin{eqnarray*}}
\def\eeann{\end{eqnarray*}}

\newcommand{\rank}{{\rm rank}}

\def\ns{\hspace{-1mm}}
\newcommand{\real}{{\mathbb{R}}}

\def\spacingset#1{\def\baselinestretch{#1}\small\normalsize}
\spacingset{1.35}

\newtheorem{lemma}{Lemma}
\newtheorem{theorem}{Theorem}
\newtheorem{remark}{Remark}
\newtheorem{corollary}{Corollary}

\def\be{\begin{equation}}
\def\ee{\end{equation}}
\def\bea{\begin{eqnarray}}
\def\eea{\end{eqnarray}}
\def\beann{\begin{eqnarray*}}
\def\eeann{\end{eqnarray*}}

\def\ns{\hspace{-1mm}}
\def\im{{\rm im}\ }
\def\proof{\noindent{\bf{\em Proof:}\ \ }}
\def\QED{\mbox{\rule[0pt]{1.5ex}{1.5ex}}}
\def\endproof{\hspace*{\fill}~\QED\par\endtrivlist\unskip}

\newcommand{\ima}{\operatorname{im}}
\newcommand{\diag}{\operatorname{diag}}
\newcommand{\defi}{\stackrel{\text{\tiny def}}{=}}

\def\gsR{{\mathcal R}^{\star}}
\def\gsV{{\mathcal V}^{\star}}

\newcommand{\complex}{{\mathbb{C}}}

\def\gD{{\cal D}}
\def\gE{{\cal E}}

\def\gJ{{\cal J}}
\def\gK{{\cal K}}
\def\gL{{\cal L}}

\def\gQ{{\cal Q}}
\def\gR{{\cal R}}
\def\gS{{\cal S}}
\def\gT{{\cal T}}
\def\gU{{\cal U}}
\def\gV{{\cal V}}

\def\gX{{\cal X}}
\def\gY{{\cal Y}}
\def\gZ{{\cal Z}}

\def\bmat{\left[ \begin{array}}
\def\emat{\end{array} \right]}

\def\bmat{\left[ \begin{array}}
\def\emat{\end{array} \right]}
\def\bsmat{\left[ \begin{smallmatrix}}
\def\esmat{\end{smallmatrix} \right]}

\def\l{{\lambda}}
\def\gU{{\cal U}}
\def\gL{{\cal L}}
\def\gR{{\cal R}}

\def\gV{{\cal V}}

\def\gS{{\cal S}}
\def\gK{{\cal K}}
\def\gT{{\cal T}}

\def\gX{{\cal X}}

\def\v{{v}}
\def\w{{w}}
\def\x{{{x}}}
\def\u{{u}}

\def\e{{e}}

\newcommand{\ba}{\begin{array}}
\newcommand{\ea}{\end{array}}

\newcommand{\spanR}{\operatorname{span}}
\def\NN{\mathbb{N}}

\begin{document}
\title{\LARGE{New geometric results in eigenstructure assignment}}

\author{Fabrizio Padula, and Augusto Ferrante, and Lorenzo Ntogramatzidis
 \thanks{Fabrizio Padula and L. Ntogramatzidis are with the Department of Mathematics and
Statistics, Curtin University, Perth,
Australia. E-mail: {\tt \{Fabrizio.Padula,L.Ntogramatzidis\}@curtin.edu.au. }}

 \thanks{Augusto Ferrante is with the Dipartimento di Ingegneria dell'Informazione, Universit\`a di Padova, 
	via Gradenigo, 6/B -- I-35131 Padova, Italy. E-mail: {\tt augusto@dei.unipd.it. }}



}


\maketitle

\vspace{-1cm}

\IEEEpeerreviewmaketitle

\begin{abstract}
The focus of this paper is the connection between two foundational areas of LTI systems theory: geometric control and eigenstructure assignment. In particular, we study the properties of the null-spaces of the reachability matrix pencil and of the Rosenbrock system matrix, which have been extensively used as two computational building blocks for the calculation of pole placing state feedback matrices and pole placing friends of output-nulling subspaces. Our objective is to show that the subspaces in the chains of kernels obtained in the construction of these feedback matrices 
interact with each other in ways that are entirely independent from the choice of eigenvalues. So far, these chains of subspaces have only been studied in the case of stationarity. In this case, it is known that these chains converge to the classic Kalman reachable subspace $\gR$ for the reachability matrix pencil and to the largest reachability subspace $\gR^\star$ in the case of the Rosenbrock matrix, respectively. Here we are interested in showing that even before stationarity has been reached, the partial chains are linked to structural properties of the system, and are therefore independent of the closed-loop eigenvalues that we wish to assign.
We further characterize these subspaces by investigating the notion of largest subspace on which it is possible to assign the closed-loop spectrum (possibly maintaining the output at zero) without resorting to non-trivial Jordan forms.
\end{abstract}


\section{Introduction}
\label{secintro}

Two fundamental frameworks that, traditionally, have been employed to study the properties of linear time invariant (LTI) dynamical systems are the so-called polynomial and the geometric approach. The area known as eigenstructure assignment, which, roughly speaking, seeks to design feedback matrices by maximizing the freedom in the assignment of the closed-loop eigenstructure, sits in between these two approaches. 

The paper \cite{Moore-76} can be considered as the initiator of this lively stream of research since, for the first time conditions were presented -- in what we will refer to as Moore's theorem \cite[Prop.~1]{Moore-76} -- outlining the freedom in the selection not only of the closed-loop eigenvalues, but also in the choice of the corresponding eigenvectors. 
Not surprisingly, generalizations of Moore's theorem have been proposed for different types of systems and with different objectives. Examples of applications include fault detection \cite{Chen-P-99}, eigenstructure assignment by static output feedback \cite{Kimura-75}, and, very importantly in the context of this paper, the construction of friends for output-nulling, reachability and stabilizability subspaces (and their duals) in geometric control theory \cite{Moore-L-78}.
 
 Loosely speaking, a pole placing feedback matrix can be obtained by firstly choosing a desired set of closed-loop eigenvalues $\l_1,\l_2,\ldots,\l_h$, and finding bases for the kernels of the $h$ matrices $[\ba{cc} \!\!A-\l_i\,I\! & B\!\!\ea]$, $i=1,\dots,h$. The first $n$ coordinates of these vectors are placed into matrices, usually called $V_i$, and the last $m$ coordinates become part of other matrices $W_i$. 
 We obtain, in this way, two matrices $V=[\ba{cccc} \!\!V_1\! &\! V_2 \!&\! \ldots\! &\! V_h\!\!\ea]$ and $W=[\ba{cccc}\!\!W_1 \!& \!W_2\! & \!\ldots\! &\! W_h\!\!\ea]$ which can be used to determine a state feedback matrix which assigns 
 the values $\l_i$ as closed-loop eigenvalues and the corresponding columns of $V$ as the closed-loop eigenvectors. More precisely, we can
 select from the columns of $V$  a basis for the Kalman reachable subspace $\gR$ (and these will become closed-loop eigenvectors), and by also selecting the corresponding columns of $W$ we obtain two matrices that we use to compute the feedback matrix (as $F=W\,V^\dagger$, where $^\dagger$ denotes the Moore-Penrose pseudo-inverse) that places the closed-loop eigenvalues at the values of the $\l_i$ that correspond to our choice of closed-loop eigenvectors. 
 This procedure can easily be adapted with some increase in the number of underlying technicalities (but not in terms of conceptual difficulty) to the defective case. This method for the construction of the pole placing feedback matrix has been recently proved to lead to much better results, in terms of robustness, precision and minimal gain, than the most widespread robust pole placement methods, including the method in \cite{Kautsky-ND-85} that MATLAB$^{\textrm{\tiny{\textregistered}}}$ uses in the routine {\tt place.m}.
 
  The procedure outlined above for pole placing feedback matrices has a direct counterpart in the case of pole placing friends of output-nulling reachability subspaces. The difference lies in the fact that the matrices $V_i$ and $W_i$ mentioned above are constructed by obtaining bases for the kernels of the Rosenbrock matrix pencil $\bsmat A-\l\,I & B \\[1mm] C & D \esmat$ in place of the reachability matrix pencil, see \cite{Moore-L-78} and \cite{NS-SICON-14}.
 
 These eigenstructure assignment techniques are important for two reasons. First, they are exhaustive: if a given set of eigenvalues can be assigned in the desired closed-loop spectrum, the feedback matrices that achieve this objective can be computed as above. Second, these techniques assign the closed-loop eigenvalues, together with the closed-loop eigenvectors. The assignment of the closed-loop eigenvectors corresponds to the possibility of shaping the response by distributing the closed-loop modes among the output components. Therefore, the possibility of assigning the closed-loop eigenstructure has been proved to be critical in the solution of several control problems, including monotonic tracking control \cite{NTSF-TAC-15}, fault detection \cite{Wahrburg-A-13} and in the context of structural decoupling  \cite{NS-SICON-18}.

These two families of pole assignment methodologies -- the one for  state feedback and the other for friends of output-nulling subspaces --  have, so far, been treated as computational techniques, because of their advantages in comparison with other methods that have been proposed in the literature. 
Indeed, only sporadically have these techniques been investigated as tools to unveil structural properties of the underlying system \cite{Aling-S-84}.
For this reason, most of the literature on these topics is restricted to the investigation of the properties of the subspaces obtained in this fashion when stationarity is obtained.
On the other hand, in areas such as tracking control, input-output decoupling and fault detection it {is} crucial to study how the column-spaces of the matrices $V_i$ and $W_i$ interact for different eigenvalues, \cite{NTSF-TAC-15,NS-SICON-18,Wahrburg-A-13}.
In fact, the freedom of selecting the closed-loop eigenvectors is reflected in the freedom of adjusting the distribution of modes among the output components (including possibly uncontrollable modes), which is crucial {in} a number of fundamental control and estimation problems.

 The purpose of this paper is to show that the subspaces obtained at each step with the methods of \cite{Moore-76} and \cite{Moore-L-78} shed light into structural properties of the system. The upper coordinates of a basis matrix of the null-space of the reachability matrix pencil or of the Rosenbrock matrix pencil span a subspace that can be imagined to rotate in the state space as a function of $\lambda$. When combining two matrices $V_i$ and $V_j$ obtained from two different eigenvalues $\l_i$ and $\l_j$, the resulting span is also $\l$-dependent: its orientation in the state-space depends on $\l_i$ and $\l_j$. However, surprisingly, its dimension is $\l$-independent: unlike what was suggested in \cite{SN-AUT-10,NS-SICON-14,NTSF-TAC-15}, once the uncontrollable eigenvalues (or, respectively, the invariant zeros) are excluded,  there are no ``bad'' choices of $\l_i$ and $\l_j$  that cause a drop in dimension of $\ima [\ba{cc} \!\! V_i & V_j\!\! \ea]$.\footnote{Clearly, from a computational point of view numerical issues
 may arise: two very close values of $\l_i$ and $\l_j$ yield subspaces $\ima V_i$ and $\ima V_j$ which are almost coincident, thus giving rise to an ill-conditioned basis for the corresponding reachable (or, respectively, reachability) space,  with the consequent numerical fragility in the computation of the state-feedback matrix or friend.}
  In this paper we show that, given $h$ distinct values $\l_1,\l_2,\ldots,\l_h$ (different from the uncontrollable eigenvalues associated with the pair $(A,B)$), the dimension of the space spanned by the columns of $[\ba{cccc} \!\!V_1 \!&\! V_2 \!&\! \ldots \!&\! V_h\!\!\ea]$ is equal to the rank of $[\ba{cccc} \!\!B &A\,B & \ldots & \!A^{h-1}B\!\!\ea]$, and it is therefore $\l$-independent for any $h$. 
  An important implication of our result is that the column-space of $[\ba{cccc} \!\!V_1 \!& \!V_2\! &\! \ldots \!& \!V_h\!\!\ea]$ is the largest closed-loop eigenspace that can be assigned with eigenvalues $\l_1,\l_2,\ldots,\l_h$ with a diagonalizable closed-loop map.
    We prove that a similar result holds true for reachability output-nulling subspaces. This is indeed a fundamental problem in control theory since, as is well known, Jordan forms are numerically ill-conditioned, and in 
 problems such as dead-beat feedback control (or filtering) a compromise realistically needs to be accepted between multiplicity of closed-loop eigenvalues which are exactly at 0 and absence of Jordan forms in the closed-loop.
%
%
We further characterize the column-space of the matrix $[\ba{cccc}\!\!V_1 \!&\! V_2 \!&\! \ldots \!&\! V_h\!\!\ea]$ by showing that its reachability subspace (i.e., the states of this subspace that can be reached from the origin with trajectories entirely contained in it)
 is a structural invariant that depends exclusively on $h$.
 Indeed, while as aforementioned the orientation of the column-space of $[\ba{cccc} \!\!V_1 \!&\! V_2 \!&\! \ldots \!&\! V_h\!\!\ea]$ depends on $\l_1,\l_2,\ldots,\l_h$ used to compute $V_1,V_2,\ldots,V_h$, the reachability on this space is independent from $\l_1,\l_2,\ldots,\l_h$. 
 In other words, all maximal output nulling subspaces obtained by assigning $h$ distinct eigenvalues share a common part which is indeed a reachability subspace.
 What happens when the rank of $[\ba{cccc} \!\!V_1\! &\! V_2 \!&\! \ldots\! &\! V_h\!\!\ea]$ reaches its stationarity now becomes obvious: in this case, we recover the classic Moore-Laub's result, i.e., the column-space of $[\ba{cccc} \!\!V_1\! &\! V_2 \!&\! \ldots\! &\! V_h\!\!\ea]$ coincides with its reachability subspace $\gR^\star$, and it is therefore independent from $\l_1,\l_2,\ldots,\l_h$.

{\bf Notation:} Throughout this paper, 
the image, kernel and Moore-Penrose pseudo-inverse of matrix $A$ are denoted by $\ima\,A$,  $\ker\,A$ and $A^\dagger$, respectively. 
Given a linear map $A: \gX \longrightarrow \gY$ and a subspace $\gS$ of $\gY$, we define 
$A^{-1}\,\gS=\{x \in \gX\,|\,A\,x \in \gS\}$. 
If $\gX=\gY$ and $\gJ$ is $A$-invariant, the eigenvalues of $A$ restricted to $\gJ$ are denoted by $\sigma\,(A\,| \gJ)$. If $\gJ_1$ and $\gJ_2$ are $A$-invariant and $\gJ_1\,{\subseteq}\,\gJ_2$, the map induced by $A$ on the quotient space $\gJ_2 / \gJ_1$ is denoted by $A\,| {\gJ_2}/{\gJ_1}$, and its spectrum is  $\sigma\,(A\,| {\gJ_2}/{\gJ_1})$.
The symbol $\oplus$ stands for the direct sum of subspaces. 
Given a map $A: \gX \longrightarrow \gX$ and a subspace $\gS$ of $\gX$, $\langle A\,|\, \gS \rangle$ is the smallest $A$-invariant subspace of $\gX$ containing $\gS$ and $\langle \gS\,|\, A \rangle$ is the largest $A$-invariant subspace contained in $\gS$. Given a vector $v\in \complex^n$, we use the symbol $\bar{v}\in \complex^n$ to denote the complex conjugate of $v$.



\section{Geometric preliminaries}
\label{PF}
Consider a quadruple $(A,B,C,D)$ associated with the non-strictly proper state-space (continuous or discrete-time) system
\beann
\Sigma:\;
\left\{\begin{array}{rcl}
\gD{\x}(t) \ns&\ns = \ns&\ns A\,\x(t)+B\,\u(t), \qquad x(0)=x_0\in \real^n \\
y(t) \ns&\ns = \ns&\ns C\,\x(t)+D\,\u(t),\end{array}\right.
\eeann
where 
$\gD$ denotes either the time derivative in the continuous time
or the unit time shift
 in the discrete time. For all $t\in \real$ in the continuous time and for all $t\in \mathbb{Z}$ in the discrete time, the vector $x(t)\in \gX=\real^n$ denotes the state, $u(t)\in \gU=\real^m$ is the  input and  $y(t)\in \gY=\real^p$ is the output of $\Sigma$.
 Thus, $A\in \real^{n\times n}$, $B\in \real^{n\times m}$, $C\in \real^{p\times n}$, and $D\in \real^{p\times m}$.

We denote by $\gR$ the reachable subspace of the pair $(A,B)$, which is the smallest $A$-invariant subspace containing $\ima B$, i.e., $\gR=\langle A\,|\,\ima B\rangle$. We denote by $\gQ$ the unobservable subspace of the pair $(C,A)$, which is the largest $A$-invariant subspace contained in the null-space of $C$, i.e., $\gQ=\langle \ker C\,|\,A\rangle$.
A subspace $\gV$ is said to be $(A,B)$-controlled invariant if,
for any $\x_0\in \gV$, there exists an input $\u(\cdot)$ such that the state trajectory $x(\cdot)$ remains on $\gV$ or, equivalently, if $A\,\gV\subseteq \gV+\ima B$, or, equivalently, if there exists a feedback matrix $F$ such that $(A+B\,F)\,\gV\subseteq \gV$.
%
%
A subspace $\gV$ is said to be output nulling if,
for any $\x_0\in \gV$, there exists $\u(\cdot)$ such that $x(\cdot)$ lies in $\gV$ and the output remains at zero or, equivalently, if $\bsmat A \\[1mm]
C\esmat \,\gV\subseteq (\gV\oplus 0_{\scriptscriptstyle \scriptscriptstyle \gY})+\ima \bsmat B\\[1mm] D \esmat$, or equivalently, if there exists a feedback matrix $F$ such that
$(A+B\,F)\,\gV\subseteq \gV \subseteq \ker (C+D\,F)$. Thus, the input that keeps $x(\cdot)$ on $\gV$ and the output at zero can be expressed as $\u(t)=F\,\x(t)$. 
  In this case, we say that $F$ is a {\em friend} of $\gV$. 
  We denote by $\mathfrak{F}(\gV)$ the set of friends of $\gV$. 

We denote by $\gV^\star_{\gE}$ the supremal output-nulling subspace contained in a subspace $\gE$ of $\gX$, which represents the set of $x_0\in \gX$ for which $u(\cdot)$ exists that maintains the state confined in $\gE$ and the output at zero {for all $t\geq 0$}. The subspace $\gV^\star_{\gE}$ can be obtained as the limit of the sequence $(\gV_i)_{i\,\in\, \mathbb{N}}$ defined by the recursion
\beann
\left\{\begin{array}{rclccccc}
 \!\! \gV_0 \ns&\ns = \ns&\ns \gE \\
 \!\! \gV_{i+1} \ns&\ns = \ns&\ns \bmat{c} \!\! A \!\! \\[-1mm] \!\! C \!\!  \emat^{-1}\left((\gV_i\oplus 0_{\scriptscriptstyle \scriptscriptstyle \gY})+\ima \bmat{c}  \!\! B \!\!  \\[-1mm] \!\!  D \!\!  \emat\right)\cap \gE, \qquad i\in \NN\setminus \{0\}.
\end{array}\right.
\eeann
The sequence $(\gV_i)_{i\,\in\, \mathbb{N}}$ is monotonically non-increasing and converges to $\gV^\star_{\gE}$ in at most $n-1$ steps, i.e., 
$\gV_0\supset \gV_1\supset \ldots \supset \gV_h = \gV_{h+1}=\ldots$ implies $\gV^\star_{\gE}=\gV_h$, with $h \le n-1$.
We denote $\gV^\star=\gV^\star_{\gX}$.

Given an output nulling subspace $\gV$, the reachability subspace $\gR_{\gV}$ on $\gV$ is the set of points that can be reached from the origin by means of input functions that keep the trajectory on $\gV$ and the output at zero. Given a friend $F$ of $\gV$, we can determine $\gR_{\gV}$ as
$\gR_{\gV}=\langle A+B\,F\,|\,\gV\cap B\,\ker D\rangle$.
The eigenstructures of  $A+B\,F\,|\,\gR_{\gV}$ and $A+B\,F\,|\,\frac{\gV+\gR}{\gV}$ are freely assignable with  $F\in \mathfrak{F}(\gV)$, whereas the spectra of $A+B\,F\,|\,\frac{\gV}{\gR_{\gV}}$ and $A+B\,F\,|\,\frac{\gX}{\gV+\gR}$ are fixed for all $F\in \mathfrak{F}(\gV)$.


An output nulling subspace $\gV$ for which a friend $F$ exists such that the spectrum of $A+B\,F\,|\,\gV$ is arbitrary is called a {\em reachability output nulling subspace}. The supremal reachability output nulling subspace is denoted by $\gR^\star$, and it coincides with the output nulling reachability subspace on $\gV^\star$, i.e., $\gR^\star = \gR_{\gV^\star}$. This subspace can be interpreted as the set of all initial states that are reachable from $0_{\scriptscriptstyle \gX}$ by inputs that keep the output at zero. The eigenstructure of $A+B\,F\,|\,\frac{\gV^\star}{\gR^\star}$ is the {\em invariant zero structure} of $\Sigma$. The  eigenvalues of $A+B\,F\,|\,\frac{\gV^\star}{\gR^\star}$ are the {\em invariant zeros} of the system: we denote by  $\gZ$ the set of invariant zeros. 

Most of the results on conditioned invariance are introduced by duality. The dual of a quadruple $(A,B,C,D)$ is $(A^\top,C^\top,B^\top,D^\top)$.
A subspace $\gS$ is $(C,A)$-conditioned invariant  if  $A\,(\gS\cap \ker C)\subseteq \gS$, or, equivalently, if a matrix $G$ exists such that $(A+G\,C)\,\gS\subseteq \gS$. 
A subspace $\gL$ is $(C,A)$-conditioned invariant if and only if $\gL^\perp$ is $(A^\top,C^\top)$-controlled invariant.
A subspace $\gS$ is said to be input containing  if $[\begin{array}{cc} \!\! A\!  & \! B \! \! \end{array}]\left((\gS \oplus \gU) \cap \ker [\begin{array}{cc}\! \!  C\!  & \! D\! \!  \end{array}]\right)\subseteq \gS$, or, equivalently, if there exists $G$ such that $\bsmat A+G\,C \\[1mm] B+G\,D \esmat \,(\gS\oplus \gU) \subseteq \gS$. 
A subspace $\gL$ is $(A,B,C,D)$-input containing if and only if $\gL^\perp$ is $(A^\top,C^\top,B^\top,D^\top)$-output nulling.
 The dual of $\gV^\star$ is the infimal input containing subspace $\gS^\star$, which is the limit of the sequence $(\gS_i)_{i \, \in \, \mathbb{N}}$ where
\beann
\left\{\begin{array}{rclccccc}
\gS_0 \ns&\ns = \ns&\ns 0_{\scriptscriptstyle \scriptscriptstyle \gX} \\
\gS_{i+1} \ns&\ns = \ns&\ns  [\begin{array}{cc}  A &  B  \end{array}]\left(({\gS}_i \oplus \gU)\cap\ker [\begin{array}{cc}   C  &  D   \end{array}]\right), \qquad i\in \mathbb{N}\setminus \{0\}.
\end{array}\right.
\eeann
This sequence is non-decreasing and converges to $\gS^\star$ in at most $n-1$ steps, i.e., 
$\gS_0\subset \gS_1\subset \ldots \subset \gS_h = \gS_{h+1}=\ldots$ implies $\gS^\star=\gS_h$, with $h \le n-1$.
By construction, given $x_j\in\gS_{j}$,  there exist $u_0,\ldots,u_{j-1}\in \gU$ such that 
\beann
x_1 \ns&\ns = \ns&\ns B\,u_0,\quad  D\,u_0=0,\\
x_2\ns&\ns = \ns&\ns(A+B\,F)\,x_1+B\,u_1, \quad  (C+D\,F)\,x_1+D\,u_1=0,\\
&\vdots&\\
x_j\ns&\ns = \ns&\ns(A+B\,F)\,x_{j-1}+B\,u_{j-1}, \quad  (C+D\,F)\,x_{j-1}+D\,u_{j-1}=0.
\eeann
In the discrete-time this means that 
each point in $\gS_j$ can be reached in at most $j$ iterations with an output that is zero up to the instant $j-1$. A counterpart of this insight for the continuous-time case can be achieved by resorting, for example, to the distributional setting of \cite{Trentelman-SH-01}.

\section{Eigenstructure assignment preliminaries}
\label{sec:eap}
%
%
The entire framework of eigenstructure assignment hinges on two fundamental results. The first one is Moore's Theorem, \cite[Prop.~1]{Moore-76}, which we recall here. 
We define the reachability matrix pencil as
\[
S_{\!\l}=[\ba{cc}A-\l\,I   &   B  \ea],
\]
and, for each $\l_i \in \complex$, we denote by $\bsmat V_{i} \\[1mm] W_{i} \esmat$ a basis matrix for the null-space of $S_{\!\l_i}$, partitioned conformably.

\begin{theorem}\cite[Prop.~1]{Moore-76}.
\label{moore}
Let $\l_1,\ldots,\l_n$ be a self-conjugate set of distinct complex numbers. There exists $F\in \real^{n \times m}$ such that $(A+B\,F)\,v_i=\l_i\,v_i$ for all $i\in \{1,\ldots,n\}$ if and only if 
\begin{enumerate}
\item $v_1,\ldots,v_n$ are linearly independent in $\complex^n$;
\item for all $i,j\in \{1,\ldots,n\}$, $v_i=\bar{v}_j$ whenever $\l_i=\bar{\l}_j$;
\item for all $i\in \{1,\ldots,n\}$, $v_i\in \ima V_{i}$.
\end{enumerate}
\end{theorem}

The importance of this result lies in the fact that conditions (1-3) guarantee the existence of a feedback matrix $F$ such that the vectors $v_1,\ldots,v_n$ are closed-loop eigenvectors of $A+B\,F$ with corresponding eigenvalues $\l_1,\ldots,\l_n$, \cite{Schmid-PN-14}.

The proof of this result in \cite{Moore-76} shows that, 
whenever a pair of eigenvalues $\l_i$ and $\l_{i+1}$ are complex conjugate, the corresponding matrices $V_{i}$ and $V_{i+1}$ can be selected to be complex conjugate as well, and defining $\hat{V}_{i}=\mathfrak{Re} \{V_{i}\}$ and 
$\hat{V}_{i+1}=\mathfrak{Im} \{V_{i}\}$ we can apply the constructive procedure in the proof of Theorem~\ref{moore} to build the feedback matrix $F$ with $\hat{V}_{i}$ and $\hat{V}_{i+1}$ in place of $V_{i}$ and $V_{i+1}$: in this way we can guarantee that the computed friend is a real matrix.  The same method, {\em mutatis mutandis}, can be used throughout this paper. 
Considering $\l_1,\l_2,\ldots,\l_h$ disjoint from the uncontrollable eigenvalues of $(A,B)$,
it is easy to see that:
\begin{itemize}
\item the rank of $[\ba{cccc}   V _{1}  &  V_{2}  &  \ldots  &  V_{h}   \ea]$ increases monotonically with $h$, and when it becomes stationary, say for $h=\rho_1$, the column-space of $[\ba{cccc}   V_{1}  &  V_{2}  &  \ldots  &  V_{\rho_1}   \ea]$ coincides with $\gR$;
\item consider a spanning matrix $[\ba{cccc}   V _{1}  &  V_{2}  &  \ldots  &  V_{h}   \ea]$, with $h\ge \rho_1$, for $\gR$. Let $r=\dim \gR$, and let $\{\v_1,\ldots,\v_r\}$ be a set of columns extracted from $[\begin{array}{ccccc}   V _{1}  &  V_{2}  &  \ldots  &  V_{h}   \end{array}]$ to form a basis for $\gR$, and let $\{\w_1,\ldots,\w_r\}$ denote the corresponding columns of $[\begin{array}{ccccc}   W_1  &  W_2  &  \ldots  &  W_h  \end{array}]$.
If $\v_k$ is a column of $V_j$, let us denote by $\mu_k$ the eigenvalue $\lambda_j$. 
Then, the matrix
\bea
\label{ultima}
F=[\begin{array}{ccccc}   \w_1  &  \w_2  &  \ldots  &  \w_r  \end{array}]\,[\begin{array}{ccccc}   \v_1  &  \v_2  &  \ldots  &  \v_r  \end{array}]^\dagger
\eea
is such that $\sigma(A+B\,F\,|\,\gR)$ is equal to the multi-set\footnote{Notice that we are not excluding the case where two different vectors $v_{k_1}$ and $v_{k_2}$ are extracted from the columns of the same matrix $V_j$. In this case, $\mu_{k_1}=\mu_{k_2}=\l_j$ and
the closed-loop eigenspace associated with the eigenvalue $\l_j$ has dimension greater than 1.}
 $\{\mu_1,\ldots,\mu_r\}$, and the closed-loop eigenvectors are $\{v_1,\ldots,v_r\}$, \cite[Cor.~1]{NP-17};
 \item when $h<\rho_1$, the column-space of $[\begin{array}{ccccc}    V _{1}  &  V_{2}  &  \ldots  &   V_{h}   \end{array}]$ (which is contained in $\gR$) is $(A,B)$-controlled invariant, and a friend can be computed as in the previous point, which assigns the eigenstructure of $A+B\,F$ restricted to $\ima [\begin{array}{ccccc}    V _{1}  &  V_{2}  &  \ldots  &  V_{h}    \end{array}]$;
 \item let $\gV$ be an $(A,B)$-controlled invariant. Let $F$ be such that $(A+B\,F)\,\gV\subseteq \gV$ and $\sigma(A+B\,F\,|\,\gV)=\{\l_1,\ldots,\l_h\}$ with $A+B\,F\,|\,\gV$ diagonalizable. Then, $F$ can always be constructed as in (\ref{ultima}) with $r=\dim \gV$, and where $\bsmat v_i\\[1mm] w_i \esmat$ is a vector of the kernel of $S_{ \mu_i}$.
This consideration shows that the construction of $F$ using null-spaces of the reachability matrix pencil is exhaustive \cite{Schmid-PN-14}.
\end{itemize}

A parallel theory is the one developed in \cite{Moore-L-78} for the supremal reachability subspace $\gR^\star$. In this case, instead of considering the reachability matrix pencil $S_{ \l}$, we consider the system matrix pencil (known as the {\em Rosenbrock matrix} \cite{Rosenbrock-70})\footnote{We recall that the invariant zeros are also characterized as the values of $\l$ for which the matrix $P_{ \l}$ loses rank with respect to its normal rank.}
\[
P_{ \l}=\bmat{cc} A-\l\,I & B \\ C & D \emat.
\]
Denoting by $\bsmat V_{i} \\[1mm] W_{i} \esmat$ a basis matrix for $\ker P_{\l_i}$ for each $\l_i \in \real$, considering a set  $\{\l_1,\ldots,\l_h\}$ disjoint from the invariant zeros, the column-space of $ [\begin{array}{ccccc}    V _{1}  &  V_{2}  &  \ldots  &  V_{h}   \end{array}]$ is output-nulling, and for a sufficiently large $h$, say $h=\rho_2$, such column-space is $\gR^\star$. A friend of $\gR^\star$, which assigns the desired closed-loop spectrum of $A+B\,F$ restricted to $\ima  [\begin{array}{ccccc}   V _{1}  &  V_{2}  &  \ldots  &  V_{h}   \end{array}]$, and guarantees the output-nulling condition $(C+D\,F)\,[\begin{array}{ccccc}   V _{1}  &  V_{2}  &  \ldots  &  V_{h}   \end{array}]=0$ can be obtained as in the previous case. When $h<\rho_2$, the column-space of $[\begin{array}{ccccc}    V _{1}  &  V_{2}  &  \ldots  &  V_{h}   \end{array}]$ (which is contained in $\gR^\star$) is an output-nulling subspace. Again, the construction of friends by extraction of vectors from the null-spaces of the system matrix pencil is exhaustive, see \cite{NS-SICON-14}.


{As mentioned in the Introduction, these eigenstructure assignment techniques are exhaustive: 
all the feedback matrices that assign a desired closed-loop spectrum can be computed as above. Moreover, assigning the associated eigenvectors corresponds to arbitrarily shaping the response by distributing the closed-loop modes among the output components.
This advantage has been explored in the context of monotonic tracking control \cite{NTSF-TAC-15} and in the state-to-output decoupling  \cite{NS-SICON-18}.}


Importantly, the exhaustiveness combined with the freedom of assigning the closed-loop eigenvectors leads to a synthesis of a feedback matrix that maximizes robustness, see e.g. \cite{Schmid-PN-14} and \cite{NS-SICON-14} for the case of pole placement and for the determination of a friend of output-nulling subspaces, respectively.

\section{Main results}

The entire framework of eigenstructure assignment for LTI systems hinges on the computation of the kernels of polynomial matrices evaluated at specific values $\l_i$ of the indeterminate $\l$. 
We now show that the null-spaces upon which these algorithms are built have important invariance properties that display (and are linked to) structural invariants of the system. 
The key questions that arise, and which have never been addressed satisfactorily, are the following:

\begin{itemize}
\item how is the dimension of the controlled invariant (resp. output nulling) subspace\\ $\ima [\ba{cccc}   V _{1}  &  V_{2}  &  \ldots  &  V_{h}   \ea]$ related to the values $\l_i$ that we use to compute a null-space of $S_{\l}$ (resp. $P_{\l}$)? In other words, are there ``good'' or ``bad'' choices of $\l_i$ that
affect the ability to extract a complete set of linearly independent vectors?\footnote{Intuition suggests that given $\l_1$, only for a zero (Lebesgue) measure set of values $\l_2$ will the matrix $[\ba{cc}    V_1 &  V_2  \ea]$ lose rank with respect to its normal rank. The possibility of the presence of these ``coupling pathologies'' could not be excluded using the frameworks of \cite{Schmid-PN-14} or \cite{NS-SICON-14}.  }

\item do these subspaces share a common part which is independent from $\l_i$? Is this part a structural invariant of the system?
\item to what extent does the convergence of $\ima [\ba{cccc}    V _{1}   &  V_{2}   &  \ldots  &   V_{h}    \ea]$ in $h$ to the subspaces $\gR$ (resp. $\gR^\star$) depend on the values of $\l_i$ that we use to the compute the null-spaces of $S_{\l}$ (resp. $P_{\l}$)? %
\end{itemize}

The next theorem is the first important result of this paper: it shows that for every $h$ (not necessarily greater or equal than $\rho_1$), the dimension of 
$\ima [\ba{cccc}   V _{1}  &  V_{2}  &  \ldots  &  V_{h}   \ea]$, where each $V_i$ is extracted from a basis matrix of $\ker S_{\l_i}$, is entirely independent from $\l_i$; such dimension is characterized in a system-theoretic sense, and can be interpreted as the dimension of the reachable subspace in $h$ steps.

\begin{theorem}
	\label{th1}
Let $\l_1,\ldots,\l_h$ be distinct, and disjoint from the uncontrollable eigenvalues of $(A,B)$. Let $\bsmat V_i\\[1mm]W_i\esmat$ be basis matrices for $\ker [\ba{cc}   A-\l_i\,I  &  B  \ea]$ for all $i\in\{1,\ldots,h\}$. We have
\be
\rank[\ba{ccccc}  V_1  &  V_2  &  \cdots  &  V_h  \ea]=
\rank[\ba{cccc}    B  &  A\,B  &  \cdots  &  A^{h-1}\,B   \ea].
\label{eq:ranghi}
\ee
\end{theorem}
\proof
First, observe that 
from the rank-nullity theorem and the fact that $\l_i$ is not uncontrollable, the rank of $\bsmat V_{i} \\[1mm] W_{i} \esmat$ is equal to $m$. 
We begin our proof under the initial assumption that $\l_1,\ldots,\l_h$ are not eigenvalues of $A$.
We also assume that $\rank B=m$. Both assumptions will be removed in the second part of the proof. From $\rank B=m$, 
we immediately see that $W_{i}$ is $m \times m$ and invertible: indeed, let $\omega$ be a vector of the null-space of $W_{i}$. We can post-multiply $(A-\l_i\,I)\,V_i+B\,W_i=0$ by $\omega$, and we obtain $(A-\l_i\,I)\,V_i\,\omega=0$. Since $\l_i \notin \sigma(A)$, then we must have $V_i\,\omega=0$, which implies that $\omega\in \ker \bsmat V_{i} \\[1mm] W_{i} \esmat$. Since $\bsmat V_{i} \\[1mm] W_{i} \esmat$ is full column-rank, we conclude that $\omega=0$.

Let us define $\bar{V}_i = V_i\,W_i^{-1}$. Clearly, $\bsmat V_{i} \\[1mm] W_{i} \esmat$ is a basis matrix for $\ker [\ba{cc}   A-\l_i\,I  &  B   \ea]$ if and only if such is $\bsmat \bar{V}_{i} \\[1mm] I \esmat$. Moreover, from 
{\small
\beann
&&\hspace{-0.1cm} \rank [\ba{cccc}   \bar{V}_1  &  \bar{V}_2  & \ldots  &  \bar{V}_{h}   \ea]= \rank\!\! \left(\!\!\![\ba{cccc}    {V}_1    & {V}_2    &\ldots    & {V}_{h}    \ea] \!\!\!
\bmat{cccc} \!\!W_1^{\scriptscriptstyle -1} \!\!&\!\! 0 \!&\! \ldots\! & \!0\!\! \\ \!\!0 \!&\! \!\!W_2^{\scriptscriptstyle -1}\! \!\!&\! \ldots \!&\! 0 \!\!\\ \!\!\vdots \!&\!\vdots\! &\! \ddots \!&\! \vdots\!\!\\ \!\!0\! & \!0\! & \!\ldots \!& \!\!W_{h}^{\scriptscriptstyle -1}\!\!\!\emat\right),
\eeann
}
 we can assume without any loss of generality that the basis matrix of $\ker [\ba{cc} A-\l_i\,I & B \ea]$ is in the form $\bsmat V_i \\ I\esmat$ and prove the statement in this case.
From $(A-\l_i\,I)\,V_i+B=0$ we find that, by defining $\mu_i=-\l_i$ and $M_i=A+\mu_i\,I$, and remembering that $\l_i\notin \sigma(A)$, we have $V_i=-M_i^{-1}\,B$. Since the matrices $M_i$ commute with each other, {we can write, 
\beann
 [\ba{cccc} \!\! V_1 \!&\! V_2 \!& \!\ldots\! &\! V_h\!\! \ea]\ns&\ns = \ns&\ns 
 L \left[\ba{cccc}
\!\! \left({\displaystyle \prod_{\small \substack{i=1\\ i\neq 1}}^h} M_i\right) B &  \left({\displaystyle \prod_{\small \substack{i=1\\ i\neq 2}}^h} M_i\right) B & \dots & \left({\displaystyle \prod_{\small \substack{i=1\\ i\neq h}}^h} M_i\right)  B \!\! \ea\right]\\
\ns&\ns = \ns&\ns  L\,[\ba{cccc}\!\!B \!& \!A\,B \!& \!\cdots \!& \!A^{h-1}\,B\!\!\ea]\,R,
\eeann
where $L\defi -\prod_{i=1}^h M_i^{-1}$ and 
$R$ is a $h \times h$ block matrix whose blocks are all scalar matrices of size $n$.
The block in row $i$ and column $k$ is given by:
\[
R_{i,k}=\left\{\begin{array}{ll}{\displaystyle \sum_{\small \substack{1\leq \ell_1<\ell_2<\dots<\ell_{h-i}\leq h\\ \ell_{j}\neq k}}} \mu_{\ell_1}\mu_{\ell_2}\ldots\mu_{\ell_{h-i}} I_m,&{\rm\ if\ } i<h\\
I_m,&{\rm\ if\ } i=h. \end{array}\right.
\]
Now we can extract the identities from the matrix $R$: by resorting to the Kronecker product, we can write
$$
R=\bar{R}\otimes I_m,
$$
with $\bar{R}$ being an $h\times h$ matrix whose entries are
\[
\bar{R}_{i,k}=\left\{\begin{array}{ll}{\displaystyle \sum_{\small \substack{1\leq \ell_1<\ell_2<\dots<\ell_{h-i}\leq h\\ \ell_{j}\neq k}}} \mu_{\ell_1}\mu_{\ell_2}\ldots\mu_{\ell_{h-i}} &{\rm\ if\ } i<h\\
1&{\rm\ if\ } i=h. \end{array}\right.
\]
Therefore, using the properties of the Kronecker product and the assumption that $\mu_j\neq \mu_i$ for $j\neq i$, we find
\[
\det R=(\det \bar{R})^m=\left(\prod_{1\leq i<j\leq h} (\mu_j-\mu_i)\right)^m \neq 0.
\]
}

We now relax the assumption $\rank B=m$. We observe that the matrix $\bsmat V_{i} \\[1mm] W_{i} \esmat$ is a basis matrix of $\ker [\ba{cc}   A-\l_i\,I  &  B  \ea]$ if and only if the matrix  $\bsmat V_{i} \\[1mm] T^{-1}\,W_{i} \esmat$ is a basis matrix for $\ker [\ba{cc}   A-\l_i\,I  &  B\,T   \ea]$ for every invertible $m \times m$ matrix $T$. Thus, no generality is lost by assuming that $B$ is either full column-rank (and in that case we have already proved the statement), or
it has the form $B=[\ba{cc}   \bar{B}  & 0   \ea]$, which is the case that we now consider. Let $\bsmat \bar{V}_i \\[1mm] \bar{W}_i \\[1mm] \bar{Z}_i \esmat$ be a basis matrix for $\ker [\ba{ccc}   A-\l_i\,I  &  \bar{B}  &  0   \ea]$. Thus
\beann
\bmat{c} \bar{V}_i \\ \bar{W}_i \\ \bar{Z}_i\emat=\bmat{cc} V_i & 0 \\ W_i & 0 \\ 0 & I \emat\,U_i,
\eeann
where $U_i=\bsmat U_{i,1}\\[1mm] U_{i,2} \esmat$ is an invertible matrix partitioned conformably (and therefore in particular $U_{i,1}$ is right-invertible), and where $\bsmat V_i \\[1mm] W_i \esmat$ is any basis matrix of $\ker [\ba{cc}   A-\l_i\,I &  B   \ea]$. We obtain the identity
\beann 
\underbrace{[\ba{cccc}      \bar{V}_1    &    \bar{V}_2    &    \ldots    &    \bar{V}_h       \ea]}_{\bar{V}}=
\underbrace{[\ba{cccc}       {V}_1    &    {V}_2    &    \ldots    &    {V}_h       \ea]}_{{V}} \underbrace{\bmat{cccc}        U_{1,1}    &    0 &    \ldots    &    0       \\[-1.5mm]       0    &    U_{2,1}    &    \ldots    &    0       \\[-1.5mm]       \vdots    &    \vdots    &    \ddots    & \vdots      \\[-1.5mm]
     0    &    0    &    \ldots & U_{h,1}       
\emat}_{U}.
\eeann
Since every $U_{k,1}$ is right-invertible, such is also $U$, and therefore we have $\bar{V}=V\,U$ and $V=\bar{V}\,U^{-R}$ (where $U^{-R}$ denotes the right inverse of $U$), which show that $V$ and $\bar{V}$ have the same rank. It follows that if $B$ is not full column-rank, we can recast the problem into one where we have the full column-rank matrix $\bar{B}$ in place of $B$. Thus, it is not restrictive to assume $\rank B=m$.

We now show that assuming that $\l_i$ are distinct from the eigenvalues of $A$ does not cause any loss of generality. We observe that $\bsmat V_i \\[1mm] W_i\esmat$ is a basis matrix of $\ker [\ba{cc}   A-\l_i\,I  &  B   \ea]$ if and only if $\bsmat V_i \\[1mm] W_i-K\,V_i\esmat$ is a basis matrix of $\ker [\ba{cc}   A+B\,K-\l_i\,I  &  B   \ea]$. Indeed,
\[
[\ba{cc}   A+B\,K-\l_i\,I  &  B \ea]=[\ba{cc}   A-\l_i\,I  &  B  \ea]\bmat{cc} I & 0 \\[-1.5mm] K & I \emat.
\]
Hence, given the distinct self-conjugate set $\{\l_1,\ldots,\l_h\}$  disjoint from the uncontrollable eigenvalues of $(A,B)$, we can always determine $K$ such that $\l_1,\ldots,\l_h$ are not eigenvalues of $A+B\,K$. From \\[-1.1cm]
{\small
\beann
&& \hspace{-0.1cm} [\ba{cccc}
B & (A+BK) B & \ldots & (A+BK)^{h-1} B \ea]\\
&&\hspace{-0.1cm} =
[\ba{cccc}
   B   &   AB  &   \ldots   &  A^{h-1} B   \ea] 
\bmat{cccccc}
  I    &    K B  &    K (A \!+\! BK) B   &   \ldots   &  K (A \!+ \!BK)^{h-1} B \!  \\[-1.5mm]
  0 &  I  &  K B  &   \ldots   & \vdots  \\[-1.5mm]
  0 & 0 & I &   \ldots   &  \vdots \\[-1.5mm]
  \vdots & \vdots & \vdots &   \ddots   &  \vdots  \\
  0 & 0 & 0 &   \ldots   &  I   \emat\!,
\eeann
}
{for any $K$ the rank of $[\ba{cccc}
  B  &  (A + BK) B  &   \ldots  &  (A + BK)^{h-1} B   \ea]$ is equal to $\rank [\ba{cccc}  
B  &  AB  &  \ldots  &  A^{h-1} B   \ea]$,} and
we can use the first part of the proof with $A+B\,K$ in place of $A$.
\endproof

Consider the reachable subspace $\gR$ of the pair $(A,B)$.
When $h\ge \rho_1$, i.e., when the rank of $[\ba{ccccc}  V_1  &  V_2  &  \cdots  &  V_h  \ea]$ becomes stationary, it is not only true that $\rank[\ba{ccccc}  V_1  &  V_2  &  \cdots  &  V_h  \ea]=\dim \gR$, but also that the identity $\ima [\ba{cccc}    B  &  A\,B  &  \cdots  &  A^{h-1}\,B   \ea]=\gR$ holds. It is well-known that the characteristic polynomial of the closed-loop matrix $A+BF$ restricted to $\gR$ can be assigned arbitrarily. One might, for instance, want to assign a single closed-loop reachable eigenvalue, i.e. $\sigma(A+B\,F\,|\,\gR)=\{\l_1\}$ (with an algebraic multiplicity equal to the dimension of $\gR$). This, however, is not always possible without using Jordan chains, or, in other words, without rendering the mapping $A+B\,F\,|\,\gR$ defective. It is therefore interesting {-- and important in contexts such as, e.g., deadbeat control --}  to ask ourselves what is the minimal cardinality of 
$\sigma(A+B\,F\,|\,\gR)$ that can be achieved with a closed-loop map $A+B\,F\,|\,\gR$ that remains diagonalizable. The next lemma sheds light onto this problem by showing an interesting correlation between the minimum number of iterations required to compute the reachable subspace of a pair $(A,B)$, and the minimum number of closed-loop reachable eigenvalues that can be assigned without the need for non-trivial Jordan forms. In other words, we are interested in characterizing the minimum  of the set
$\gT_1=\big\{q\in \NN\,|\,\exists\,F\in \real^{m \times n}:\; \text{$A+B\,F\,|\,\gR$ is diagonalizable and $\operatorname{card}\bigl(\sigma (A+B\,F\,|\,\gR)\bigr)=q$}\big\}$.

\begin{lemma}
	\label{lem:lemma2}
	Let $h=\min \{\ell\in \NN\,|\, \gR=\ima [\begin{array}{cccc}    B   &   A\,B  &  \cdots   &  A^{\ell-1} B     \end{array}]\}$. Then, $h=\min \gT_1$.
	\end{lemma}
%
%
%
%
\proof 
 Obviously $h \leq \dim \gR$. 
In view of Theorem~\ref{th1}
we have $\gR=\rank[\ba{ccccc}  V_1  &  V_2  &  \cdots  &  V_h  \ea]$ and we can always find $F$ such that $\sigma (A+B\,F\,|\,\gR)=\{\l_1,\ldots,\l_h\}$, where $\l_1,\ldots,\l_h \notin \sigma(A)$, and the closed-loop restricted to $\gR$ is diagonalizable. This shows that $h\in \gT_1$. Now we show that $h$ is also the minimum of $\gT_1$. By contradiction, assume that 
there exists $\ell\in \gT_1$ such that $\ell<h$. 
Let $F$ be a matrix such that  $A+B\,F\,|\,\gR$ is diagonalizable and $\operatorname{card}\bigl(\sigma (A+B\,F\,|\,\gR)\bigr)=\ell$. Let $\{\l_1,\ldots,\l_{\ell}\}= \sigma (A+B\,F\,|\,\gR)$. 
Since the parameterization in \cite[Prop.~2.1]{Schmid-PN-14} is complete, there exist basis matrices $\bsmat V_1 \\ W_1\esmat$, $\ldots$, $\bsmat V_{\ell} \\ W_{\ell}\esmat$ of $S_{\l_1},\ldots,S_{\l_{\ell}}$ such that there exists a selection 
of $r$ columns $v_1,\ldots,v_r$ (where $r=\dim \gR$) from the columns of  $[\ba{ccccc}  V_1  &  V_2  &  \cdots  &  V_{\ell}  \ea]$ and corresponding columns $w_1,\ldots,w_r$ from $[\ba{ccccc}  W_1  &  W_2  &  \cdots  &  W_{\ell}  \ea]$ that allows us to define $F=[\ba{ccccc}  w_1  &  w_2  &  \cdots  &  w_{r}  \ea]\,[\ba{ccccc}  v_1  &  v_2  &  \cdots  & v_{r}  \ea]^{\dagger}$, so that
$\gR=\spanR\{v_1,v_2,\ldots,v_r\}=
\ima [\ba{ccccc}  V_1  &  V_2  &  \cdots  &  V_{\ell}  \ea]$, which,
in view of Theorem~\ref{th1}, implies $\gR = \im [\begin{array}{cccc}   B  &  A\,B & \cdots  & A^{\ell-1}\,B   \end{array}]$, which contradicts the fact that $h$ is the minimum of the set $\{\ell\in \NN\,|\, \gR=\ima [\begin{array}{cccc}    B   &   A\,B  &  \cdots   &  A^{\ell-1} B     \end{array}]\}$. 
\endproof

We now parallel Theorem~\ref{th1} with a result that characterizes the dimension of the output-nulling subspaces obtained by matrices which are extracted from bases of the null-space of $P_{\l_i}$. As it will be clear in the sequel, this dimension equals the dimension of the space of states that can be reached in $h$ steps by maintaining the output at zero and from which the system can evolve with zero output.

First, we introduce two preliminary results. The first is proved in \cite[p. 170]{Trentelman-SH-01}.

\begin{lemma}
\label{basis}
Consider a change of basis matrix $T=[\begin{array}{ccc}   T_1  &  T_2  &  T_3   \end{array}]$ in $\gX$ such that the columns of $T_1$ are a basis for $\gR^\star$ and the columns of $T=[\begin{array}{ccc}     T_1   &   T_2    \end{array}]$ are a basis for $\gV^\star$. 
Consider a change of basis $\Omega=[\begin{array}{ccc}     \Omega_1   &   \Omega_2     \end{array}]$ in $\gU$ such that the columns of $\Omega_1$ are a basis for $B^{-1}\,\gV^\star \cap \ker D$.
Given any friend $F$ of $\gV^\star$ (which is also a friend of $\gR^\star$),
{\small 
\beann
\bar{A} \ns&\ns    =    \ns&\ns  T^{-1}(A  \!+\!  B F)\,T=\bmat{ccc}
   \!  \bar{A}_{\scriptscriptstyle 1,1}   &   \bar{A}_{\scriptscriptstyle 1,2}   &   \bar{A}_{\scriptscriptstyle 1,3}   \!  \\
  \!  0   &   \bar{A}_{\scriptscriptstyle 2,2}   &   \bar{A}_{\scriptscriptstyle 2,3} \!    \\ 
  \!  0   &   0   &   \bar{A}_{\scriptscriptstyle 3,3} \!    \emat\!\!,  \! \!\!  \quad \bar{B}=T^{-1} B\, \Omega=\bmat{cc} 
   \! \bar{B}_{\scriptscriptstyle 1,1}   &   \bar{B}_{\scriptscriptstyle 1,2}\!    \\
\!0   &   \bar{B}_{\scriptscriptstyle 2,2}  \!  \\
\!0   &   \bar{B}_{\scriptscriptstyle 3,2}  \!  \emat \\
\bar{C}    \ns&\ns    =    \ns&\ns  (C+D\,F)\,T=[\ba{ccc}0 & 0 & \bar{C}_{\scriptscriptstyle 3}\ea], \quad 
\bar{D}=D\,\Omega=[\ba{cc}0 & \bar{D}_{\scriptscriptstyle 2}\ea].
\eeann
}
In this basis the following facts hold:
\begin{itemize}
\item If we denote by $\gV^\star$ and $\gS_j$ the supremal output-nulling and the $j$-th term of the sequence of $\gS^\star$  of the quadruple $(A,B,C,D)$, and by $\bar{\gV}^\star$ and $\bar{\gS}_j$ the corresponding subspaces of the quadruple $(\bar{A},\bar{B},\bar{C},\bar{D})$, we have $\bar{\gV}^\star=T^{-1}\,\gV^\star$ and $\bar{\gS}^\star=T^{-1}\,\gS^\star$; 
\item the invariant zero structure is the eigenstructure of $\bar{A}_{2,2}$; 
\item the matrix $\bsmat \bar{B}_{3,2}\\[1mm] \bar{D}_2 \esmat$ is full column-rank;
\item the pair $(\bar{A}_{1,1},\bar{B}_{1,1})$ is completely reachable;
\item the supremal output-nulling subspace $\bar{\gV}^\prime$ of the subsystem $(\bar{A}_{3,3},\bar{B}_{3,2},\bar{C}_{3},\bar{D}_2)$ is $\{0\}$; 
\item the Smith form of 
$\bsmat \bar{A}_{3,3}-\l\,I & \bar{B}_{3,2}\\[1mm]
\bar{C}_{3} & \bar{D}_2\esmat
$ is $\bsmat I \\[1mm] 0 \esmat$. 
\end{itemize}
\end{lemma}

We now present the second preliminary result, see \cite[p. 690]{Commault-DDLM-86}.

\begin{lemma}
\label{intersection}
For all $i,j\in \mathbb{N}$ there holds
{\small
\beann
&&\hspace{-4mm} \gV_i  \cap\gS_j  = [\begin{array}{ccccccccc} \! A^{\scriptscriptstyle  j-1} B \! &\!   \ldots  \!&\!  A B  \! &   B  &\!  0\!  &\!  \ldots\! \!  &  0\!   \end{array}] \ker \!
\bmat{ccccccc}
D & 0 &  \!  \!  \! \ldots \!  \!  \! & 0 & 0  \!  \!\!  \\
 \!  \! C\,B & D &  \!  \! \!  \ldots \!  \!  \! & 0 & 0  \!  \! \! \\[-2mm]
\vdots & \vdots &  \! \!  \! \ddots \!  \!  \! & \vdots \! \! \! \\[-1.5mm]
 \! \!  \! C\,A^{\scriptscriptstyle j-2}\,B \! \!  \! & \!   C\,A^{\scriptscriptstyle j-3}\,B   \! & \!  \! \!  \ldots  \! \!  \! & \!  0  \! &  \! 0  \!  \! \! \\[-2mm]
\vdots & \vdots & \vdots & \!  \!\!   \ddots  \! \!  \! & \vdots \! \! \\[-1.5mm]
\!  C A^{\scriptscriptstyle i+j-2} B   \! &   C A^{\scriptscriptstyle i+j-3} B \! &\!     \ldots\!   & \! C B\!   &\!  D \!  \emat\!\!.
\eeann
}
\end{lemma}
The generalization of Theorem~\ref{th1} is as follows.

\begin{theorem}
	\label{th2}
Let $\l_1,\ldots,\l_h$ be distinct and disjoint from the invariant zeros of $\Sigma$ and let $\bsmat V_i\\[1mm]W_i\esmat$ be basis matrices for $\ker P_{\!\l_i}$ for all $i\in\{1,\ldots,h\}$. Then,
\be
\rank[\ba{ccccc}  V_1  &  V_2  &  \cdots  &  V_h  \ea]=\dim(\gV^\star\cap \gS_h).
\label{eq:ranghi2}
\ee
\end{theorem}
\proof
Consider the quadruple $(\bar{A},\bar{B},\bar{C},\bar{D})$ in Lemma~\ref{basis}. 
We use Lemma~\ref{intersection} with $i=n$ and $j=h$, so that it is guaranteed that $\gV_i=\gV^\star$.
In this basis it is straightforward to see that
$\bar{C}\,\bar{A}^k\,\bar{B} = [\ba{ccc}    0  &  \bar{C}_{3}\,\bar{A}_{3,3}^k\,\bar{B}_{3,2}  \ea]$, 
so that
{\small
\bea
 \ker \bmat{cccccccc} 
\bar{D} & 0 & 0 & \ldots & 0 \\
\bar{C}\,\bar{B} & \bar{D} & 0 & \ldots& 0 \\
\bar{C}\,\bar{A}\,\bar{B} & \bar{C}\,\bar{B} & \bar{D} & \ldots & 0  \\
\vdots & \vdots & \vdots & \ddots & \vdots \\
\bar{C}\,\bar{A}^{\scriptscriptstyle n+h-2}\,\bar{B} & \bar{C}\,\bar{A}^{\scriptscriptstyle n+h-3}\,\bar{B} &\ldots & \ldots & \bar{D}
\emat \ns&\ns = \ns&\ns \ker 
\bmat{cc|cc|c|cc}
0 & \bar{D}_{\scriptscriptstyle 2} & 0 & 0 & \ldots & 0 & 0 \\
\hline
0 &  \bar{C}_{\scriptscriptstyle 3}\,\bar{B}_{3,2} & 0 &\bar{D}_{\scriptscriptstyle 2} &  \ldots & 0 & 0 \\
\hline
0 & \bar{C}_{\scriptscriptstyle 3}\,\bar{A}_{3,3}\,\bar{B}_{\scriptscriptstyle 3,2} & 0 &  \bar{C}_{\scriptscriptstyle 3}\,\bar{B}_{\scriptscriptstyle 3,2} &  \ldots & 0 & 0 \\
\hline
\vdots & \vdots &  \vdots & \vdots &  \ddots & \vdots & \vdots \\
\hline
0 & \bar{C}_{\scriptscriptstyle 3}\,\bar{A}_{\scriptscriptstyle 3,3}^{\scriptscriptstyle n+h-2} \,\bar{B}_{\scriptscriptstyle 3,2} & 
0 & \bar{C}_{\scriptscriptstyle 3}\,\bar{A}_{\scriptscriptstyle 3,3}^{\scriptscriptstyle n+h-3} \,\bar{B}_{\scriptscriptstyle 3,2} & \ldots & 0& \bar{D}_{\scriptscriptstyle 2} \emat \nonumber \\
 \ns&\ns = \ns&\ns \ima \bmat{c|c|c|c|c|c}
\; I \;& \;0 \;&\; 0 \;&\;0\; &\; \ldots\; & \;0\; \\[-1.0mm]
0 & 0 & 0 & 0 & \ldots & K_{\scriptscriptstyle 0} \\[-0.7mm]
\hline
0 & I & 0 & 0 & \ldots & 0 \\[-1.0mm]
0 & 0 & 0 & 0 & \ldots & K_{\scriptscriptstyle 1} \\[-0.7mm]
\hline
\vdots & \vdots & \vdots & \vdots &  \ldots & \vdots \\[-0.7mm]
\hline
0 & 0 & 0 & 0 & \ldots & 0 \\[-1.0mm]
0 & 0 & 0 & 0 & \ldots & K_{\scriptscriptstyle n+h-1} \!\!
\emat\label{id}
\eea
}
where $\bsmat K_0 \\ K_1 \\[-2mm] \vdots \\ K_{n+h-1} \esmat$ is a basis of $\ker M$, with
\beann
M\defi
\bmat{ccccc}
\bar{D}_{\scriptscriptstyle 2} & 0 & 0 & \ldots & 0 \\
  \bar{C}_{\scriptscriptstyle 3}\,\bar{B}_{\scriptscriptstyle 3,2} & \bar{D}_{\scriptscriptstyle 2} & 0 & \ldots & 0 \\
 \bar{C}_{\scriptscriptstyle 3}\,\bar{A}_{\scriptscriptstyle 3,3}\,\bar{B}_{\scriptscriptstyle 3,2}  &  \bar{C}_{\scriptscriptstyle 3}\,\bar{B}_{\scriptscriptstyle 3,2} & \bar{D}_{\scriptscriptstyle 2} &  \ldots  & 0 \\
\vdots & \vdots& \vdots &  \ddots &  \vdots \\
 \bar{C}_{\scriptscriptstyle 3}\,\bar{A}_{\scriptscriptstyle 3,3}^{\scriptscriptstyle n+h-2} \,\bar{B}_{\scriptscriptstyle 3,2} & 
 \bar{C}_{\scriptscriptstyle 3}\,\bar{A}_{\scriptscriptstyle 3,3}^{\scriptscriptstyle n+h-3} \,\bar{B}_{\scriptscriptstyle 3,2} & 
\bar{C}_{\scriptscriptstyle 3}\,\bar{A}_{\scriptscriptstyle 3,3}^{\scriptscriptstyle n+h-4} \,\bar{B}_{\scriptscriptstyle 3,2} & \ldots & \bar{D}_{\scriptscriptstyle 2} \emat.
\eeann
We now show that $K_0=K_1=\ldots=K_{h-1}=0$. Suppose by contradiction that $K_0\neq 0$. 
Let
\[
\bmat{c} \omega(0) \\[-2mm]  \vdots \\ \omega(n)\emat\in \ima 
\bmat{c} K_0 \\[-2mm]  \vdots \\ K_n \emat
\]
with $\omega(0)\neq 0$. Thus, $\bar{D}_2\,\omega(0)=0$, and since $\bsmat \bar{B}_{3,2} \\[1mm] \bar{D}_2\esmat$ is injective from Lemma~\ref{basis}, it follows also that $\bar{B}_{3,2}\,\omega(0) \neq 0$. Let $\xi(1)=\bar{B}_{3,2}\,\omega(0)\neq 0$. 
We obtain a contradiction by showing that $\xi(1)\in \bar{\gV}^\prime=\{0\}$. Indeed, considering $(\bar{A}_{3,3},\bar{B}_{3,2},\bar{C}_{3},\bar{D}_2)$ as a discrete-time system
\beann
\xi(k+1)\ns&\ns =\ns&\ns \bar{A}_{3,3}\xi(k)+\bar{B}_{3,2}\,\omega(k) \\
z(k) \ns&\ns  =\ns&\ns\bar{C}_{3}\,\xi(k)+\bar{D}_2\,\omega(k)
\eeann
 we have found an input $\omega(1),\omega(2),\ldots,\omega(n)$ such that the corresponding outputs $z(1),\ldots,z(n)$ are all equal to zero:
 \beann
 z(1)  \ns&\ns  =\ns&\ns \bar{C}_{3}\,\xi(1)+\bar{D}_2\,\omega(1)
=
[\ba{cc}  \bar{C}_{\scriptscriptstyle 3}\,\bar{B}_{\scriptscriptstyle 3,2}& \bar{D}_{\scriptscriptstyle 2}\ea] \bmat{c} \omega(0) \\[-1.2mm] \omega(1) \emat=0\\[-2mm]
 z(2)  \ns&\ns  =\ns&\ns [\ba{ccc}\bar{C}_{\scriptscriptstyle 3}\,\bar{A}_{\scriptscriptstyle 3,3}\,\bar{B}_{\scriptscriptstyle 3,2}& \bar{C}_{\scriptscriptstyle 3}\,\bar{B}_{\scriptscriptstyle 3,2}& \bar{D}_{\scriptscriptstyle 2}\ea] \!\!\bmat{c} \omega(0) \\ \omega(1) \\ \omega(2) \emat=0\\[-5mm]
  \ns&\ns  \vdots \ns&\ns\\[-5mm]
   z(n)  \ns&\ns  =\ns&\ns [\ba{ccccc}  \bar{C}_{\scriptscriptstyle 3}\,\bar{A}_{\scriptscriptstyle 3,3}^{\scriptscriptstyle n-1} \bar{B}_{\scriptscriptstyle 3,2}&\bar{C}_{\scriptscriptstyle 3}\,\bar{A}_{\scriptscriptstyle 3,3}^{\scriptscriptstyle n-2} \bar{B}_{\scriptscriptstyle 3,2}& \ldots & \bar{C}_{\scriptscriptstyle 3}\,\bar{B}_{\scriptscriptstyle 3,2}& \bar{D}_{\scriptscriptstyle 2}\ea] \!\!\bmat{c} \omega(0) \\[-1.2mm] \omega(1) \\[-2mm] \vdots \\[-1.2mm] \omega(n) \emat=0.
\eeann
We obtained a contradiction. Thus, we have proved that $K_0=0$. We can repeat the same argument for $K_1$. Indeed, since now $K_0=0$, the part of $M$ that we need to consider is obtained by removing the first block of rows and columns, and what remains has the same structure of $M$. This argument ends, as observed above, with the term $K_{h-1}$. 
We have proved that $K_0=K_1=\ldots=K_{h-1}=0$.
We have also 
\beann
\bar{A}^h\,\bar{B} \ns&\ns = \ns&\ns \bmat{cc} 
\bar{A}_{1,1}^h\,\bar{B}_{1,1} & \Phi_{h-1} \\
0 & \Psi_{h-1}\\
0 & \Theta_{h-1} \emat
\eeann
for suitable matrices $\Phi_{h-1}, \Psi_{h-1}$ and $\Theta_{h-1}$, and using (\ref{id}) we obtain
{\small
\beann
&& \hspace{-1mm} [\begin{array}{ccccccccc}  \bar{A}^{\scriptscriptstyle h-1} B  &   \ldots  &  \bar{A} \bar{B}   &   \bar{B}  &  0  &   \ldots   &  0   \end{array}] \ker 
\bmat{ccccccc}
\bar{D} & 0 &       \ldots      & 0 & 0      \\
    \bar{C}\,\bar{B} & \bar{D} &       \ldots      & 0 & 0      \\
\vdots & \vdots &      \ddots      & \vdots    \\
     \bar{C}\,\bar{A}^{\scriptscriptstyle h-2}\,\bar{B}     &    \bar{C}\,\bar{A}^{\scriptscriptstyle h-3}\,\bar{B}    &      \ldots      &   0   &   0      \\
\vdots & \vdots & \vdots &      \ddots      & \vdots   \\
   \bar{C} \bar{A}^{\scriptscriptstyle n+h-2} \bar{B}    &   \bar{C} \bar{A}^{\scriptscriptstyle n+h-3} \bar{B}  &     \ldots     &  \bar{C} \bar{B}   &  \bar{D}    \emat\\
&& \hspace{-3mm} = \ima \bmat{c|c|c|c|ccc|c}
\bar{A}^{\scriptscriptstyle h-1}_{\scriptscriptstyle 1,1}\,\bar{B}_{\scriptscriptstyle 1,1}& \bar{A}^{\scriptscriptstyle h-2}_{\scriptscriptstyle 1,1}\,\bar{B}_{\scriptscriptstyle 1,1}& \ldots & \bar{B}_{\scriptscriptstyle 1,1} & 0 & \ldots & 0 & H_{\scriptscriptstyle h-1}^\prime \\
0 & 0 & \ldots & 0 & 0 & \ldots & 0 & H_{\scriptscriptstyle h-1}^{\prime \prime} \\
0 & 0 & \ldots & 0 & 0 & \ldots & 0 & 
H_{\scriptscriptstyle h-1}^{\prime \prime \prime}\emat,
\eeann
}
where
\beann
H_{\scriptscriptstyle h-1}^\prime \ns&\ns = \ns&\ns 
\Phi_{\scriptscriptstyle h-1}\,K_{\scriptscriptstyle 0}+\ldots+\Phi_{\scriptscriptstyle 0}\,K_{\scriptscriptstyle h-1}+0\,K_{\scriptscriptstyle h}+\ldots+0\,K_{\scriptscriptstyle n+h-1}=0 \\
H_{\scriptscriptstyle h-1}^{\prime \prime} \ns&\ns = \ns&\ns 
\Psi_{\scriptscriptstyle h-1}\,K_{\scriptscriptstyle 0}+\ldots+\Psi_{\scriptscriptstyle 0}\,K_{\scriptscriptstyle h-1}+0\,K_{\scriptscriptstyle h}+\ldots+0\,K_{\scriptscriptstyle n+h-1}=0 \\
H_{\scriptscriptstyle h-1}^{\prime \prime \prime} \ns&\ns = \ns&\ns 
\Theta_{\scriptscriptstyle h-1}\,K_{\scriptscriptstyle 0}+\ldots+\Theta_{\scriptscriptstyle 0}\,K_{\scriptscriptstyle h-1}+0\,K_{\scriptscriptstyle h}+\ldots+0\,K_{\scriptscriptstyle n+h-1}=0. 
\eeann
It follows that
{\small
\beann
&&\hspace{-1mm} \dim\! \left( \!\!\! [\begin{array}{ccccccccc}  \bar{A}^{\scriptscriptstyle h-1} B  &  \!   \ldots  \!    &  \bar{A} \bar{B}   &   \bar{B}  &  0  &    \!   \ldots  \!     &  0   \end{array}] \ker  \! \!
\bmat{ccccccc}
\bar{D} & 0 &    \!\!    \ldots  \!\!    & 0 & 0      \\[-1mm]
    \bar{C} \bar{B} & \bar{D} &       \!\!  \ldots  \!\!      & 0 & 0      \\[-1.2mm]
\vdots & \vdots &      \!\!  \ddots     \!\!   &    \vdots &    \vdots    \\[-1mm]
     \bar{C} \bar{A}^{\scriptscriptstyle h-2} \bar{B}     &    \bar{C} \bar{A}^{\scriptscriptstyle h-3} \bar{B}    &    \! \!     \ldots    \! \!     &   0   &   0      \\[-1.2mm]
\vdots & \vdots & \! \! \ddots\! \!  &         \vdots        & \vdots   \\[-1mm]
   \!\!  \bar{C} \bar{A}^{\scriptscriptstyle n+h-2} \bar{B}    &   \bar{C} \bar{A}^{\scriptscriptstyle n+h-3} \bar{B}    &    \!\!    \ldots   \!\!    &  \bar{C} \bar{B}   &  \bar{D} \! \!    \emat\right)\\
&&=\rank [\ba{cccc} 
 \bar{A}^{\scriptscriptstyle h-1}_{\scriptscriptstyle 1,1}\,\bar{B}_{\scriptscriptstyle 1,1}& \bar{A}^{\scriptscriptstyle h-2}_{\scriptscriptstyle 1,1}\,\bar{B}_{\scriptscriptstyle 1,1}&\ldots & \bar{B}_{\scriptscriptstyle 1,1} \ea].
\eeann
}
We recall that $\bsmat V_i \\[1mm] W_i\esmat$ denotes a basis matrix of $\ker P_{\l_i}$, so that in the given basis 
\[
P_{\l_i}=\bmat{ccc|cc}
     \bar{A}_{\scriptscriptstyle 1,1}-\l_i\,I   &   \bar{A}_{\scriptscriptstyle 1,2}   &   \bar{A}_{\scriptscriptstyle 1,3}     & \bar{B}_{1,1}   &   \bar{B}_{1,2}    \\
    0   &   \bar{A}_{\scriptscriptstyle 2,2}-\l_i\,I   &   \bar{A}_{\scriptscriptstyle 2,3}     &  0   &   \bar{B}_{\scriptscriptstyle 2,2}    \\ 
    0   &   0   &   \bar{A}_{\scriptscriptstyle 3,3}-\l_i\,I     &  0   &   \bar{B}_{\scriptscriptstyle 3,2}   \\
 \hline
 0 & 0 & \bar{C}_{\scriptscriptstyle 3}&  0 & \bar{D}_{\scriptscriptstyle 2}\emat
 \]
and $\bsmat V_i \\[1mm] W_i\esmat$ is partitioned as
$\bmat{ccc|cc} V_{1,i}^\top & V_{2,i}^\top & V_{3,i}^\top &W_{1,i}^\top & W_{2,i}^\top \emat^\top$.
We show that $V_i=\bsmat V_{1,i} \\[1mm] 0\\[1mm] 0
\esmat$ and $W_i=\bsmat W_{1,i} \\[1mm] 0
\esmat$, where the columns of $\bsmat V_{1,i} \\ W_{1,i}\esmat$ are a basis for $\ker  [\ba{ccc}     \bar{A}_{1,1}-\l_i\,I   & \ \bar{B}_{1,1}    \ea]$.
To this end, observe that from the structure of $P_{\l_i}$, the columns of $\bsmat V_{3,i}\\[1mm] W_{2,i}\esmat$ are a basis for 
$\ker \bsmat \bar{A}_{3,3}-\l_i\,I  & \bar{B}_{3,2}\\
\bar{C}_{3} & \bar{D}_2\esmat$.
From Lemma~\ref{basis}, since the Smith form of $\bsmat \bar{A}_{3,3}-\l_i\,I  & \bar{B}_{3,2}\\
\bar{C}_{3} & \bar{D}_2\esmat$ is $\bsmat I \\[1mm] 0 \esmat$, we have $\bsmat V_{3,i}\\ W_{2,i}  \esmat=0$. It follows that 
$\bmat{ccc|cc} V_{1,i}^\top & V_{2,i}^\top & 0 & W_{1,i}^\top & 0 \emat^\top$
is a basis matrix of $\ker P_{\l_i}$, which implies that
$\bmat{ccc} V_{1,i}^\top & V_{2,i}^\top &  W_{1,i}^\top \emat^\top$
is a basis matrix for 
$\ker \bsmat  \bar{A}_{1,1}-\l_i\,I &   \bar{A}_{1,2}  & \bar{B}_{1,1}  \\
 0   &   \bar{A}_{2,2}-\l_i\,I   &  0   \esmat$.
Since $\l_i$ is not an invariant zero of the system, we obtain $V_{2,i}=0$, so that 
$\bmat{cc} V_{1,i}^\top &   W_{1,i}^\top \emat^\top$
is a basis matrix for 
$\ker [\ba{ccc}     \bar{A}_{1,1}-\l_i\,I   &  \bar{B}_{1,1}    \ea]$. It follows that
{\small
\beann
&&\hspace{-2mm} \dim (\gV^\star \cap \gS_h)= \dim (\bar{\gV}^\star \cap \bar{\gS}_h)\\
&& \hspace{-2mm}=  \dim \!  \!\left(  \! \! [\begin{array}{ccccccccc}  \bar{A}^{\scriptscriptstyle h-1} B  &     \ldots    &  \bar{A} \bar{B}   &   \bar{B}  &  0  &     \ldots     &  0   \end{array}]
 \ker  \! \!
\bmat{ccccccc}
\bar{D}  &  0  &  \! \!  \ldots  \! \! &  0  &  0    \\[-1mm]
 \bar{C} \bar{B} & \bar{D} &   \! \!  \ldots \! \!  & 0 & 0   \\[-2mm]
\vdots & \vdots &   \! \!\ddots  \! \!&  \vdots   &  \vdots    \\[-1mm]
   \bar{C} \bar{A}^{\scriptscriptstyle h-2} \bar{B}  &   \bar{C} \bar{A}^{\scriptscriptstyle h-3} \bar{B}   &  \! \!  \ldots   \! \! &   0   &   0      \\[-2mm]
\vdots & \vdots &  \! \!\ddots \! \! &     \vdots   &  \vdots  \\[-1mm]
 \! \! \bar{C} \bar{A}^{\scriptscriptstyle n+h-2} \bar{B}  &  \bar{C} \bar{A}^{\scriptscriptstyle n+h-3} \bar{B}  &    \! \!\ldots \! \!&\bar{C} \bar{B}  &  \bar{D}   \! \!  \emat\right)\\
&& \hspace{-2mm}= \rank [\ba{cccc} 
\bar{A}^{\scriptscriptstyle h-1}_{\scriptscriptstyle 1,1}\,\bar{B}_{\scriptscriptstyle 1,1}& \bar{A}^{\scriptscriptstyle h-2}_{\scriptscriptstyle 1,1}\,\bar{B}_{\scriptscriptstyle 1,1}& \ldots & \bar{B}_{\scriptscriptstyle 1,1} \ea] \\
&& \hspace{-2mm}= \rank [\ba{cccc} V_{\scriptscriptstyle 1,1} & V_{\scriptscriptstyle 1,2} & \ldots & V_{\scriptscriptstyle 1,h} \ea] \\
&& \hspace{-2mm}= \rank \bmat{cccc} V_{\scriptscriptstyle 1,1} & V_{\scriptscriptstyle 1,2} & \ldots & V_{\scriptscriptstyle 1,h} \\ 0 & 0 & \ldots & 0 \\  0 & 0& \ldots & 0 
\emat = \rank [\ba{cccc} V_{\scriptscriptstyle 1} & V_{\scriptscriptstyle 2} & \ldots & V_{\scriptscriptstyle h} \ea],
\eeann
}
where the fourth equality is a consequence of Theorem~\ref{th1}.
\endproof

Theorem~\ref{th1} can be interpreted as a special case of the result of Theorem~\ref{th2}, with $C$ and $D$ empty. In fact, in that case, $\gV^\star=\gX$ and $\gS_h=\ima [\ba{cccc}  B & A\,B & \ldots &  A^{h-1}\,B \ea]$. 

We now provide a geometric interpretation of the output-nulling subspace spanned by the columns of $[\ba{cccc}  V_1 & V_2 & \ldots & V_h \ea]$ obtained with a given set of closed-loop eigenvalues $\l_1,\ldots,\l_h$.

\begin{theorem}
\label{lattice}
Let $\l_1,\ldots,\l_h$ be self-conjugate, distinct, and disjoint from the invariant zeros. {Let $\gV_{\Sigma}$ denote the set of output-nulling subspaces of $\Sigma$.}
Then:
\begin{itemize}
\item The set
\beann
&& \hspace{-0.8cm} \mathfrak{T}=\big\{\gV\in \gV_{\Sigma}\,|\,\exists\,F\in \mathfrak{F}(\gV):\; \text{$A+B\,F\,|\,\gV$ is diagonalizable and $\sigma (A+B\,F\,|\,\gV)=\{\l_1,\ldots,\l_h\}$}\big\}
\eeann
admits a maximum $\gK_h$;
\item $\gK_h=\ima [\ba{cccc}  V_1 & V_2 & \ldots & V_h \ea]$.
\end{itemize}
%
\end{theorem}
\proof
Let $\gV_1,\gV_2\in \mathfrak{T}$. Let $F_1\in \mathfrak{F}(\gV_1)$ such that $A+B\,F_1\,|\,\gV$ is diagonalizable and $\sigma (A+B\,F_1\,|\,\gV_1)=\{\l_1,\ldots,\l_h\}$. We want to show that there exists $\gV\in  \mathfrak{T}$ which contains both $\gV_1$ and $\gV_2$.
Hence, there hold
\bea
(A+B\,F_1)\,V_1\ns&\ns =\ns&\ns V_1\,X \label{e11c} \\
(C+D\,F_1)\,V_1\ns&\ns =\ns&\ns 0\label{e22c}
\eea
where $X$ is diagonalizable and is similar to the matrix $\diag \{\l_1,\ldots,\l_h\}$.

Let $v_{1,1},\ldots,v_{1,\eta_1}$ be the set of closed-loop eigenvectors of $A+B\,F_1\,|\,\gV_1$.\\
 Let $T=[\ba{ccc|c}  v_{1,1}&\ldots&v_{1,\eta_1} & \star \ea]$ be a change of basis in $\gX$. We find
\bea
T^{-1}(A+B\,F_1)\,T\,T^{-1}\,V_1\ns&\ns =\ns&\ns T^{-1} V_1\,T\,T^{-1}XT \label{e11cl} \\
(C+D\,F_1)\,T\,T^{-1}\,V_1\ns&\ns =\ns&\ns 0\label{e22cl}
\eea
i.e.,
\[
\bmat{cc} A^F_{1,1} & \star \\[-1mm] 0 & \star\emat\,\bmat{c} I \\[-1mm] 0 \emat=\bmat{c} I \\[-1mm] 0 \emat\,\operatorname{diag}\{\underbrace{\l_1,\ldots,\l_1}_{\nu_1},\ldots,\underbrace{\l_h,\ldots,\l_h}_{\nu_h}\},
\]
where $\nu_1+\ldots+\nu_h=\dim \gV_1=\eta_1$ and therefore $A^F_{1,1}$ is diagonal. From this structure, in the new basis the friend $\tilde{F}_1=F_1\,T$ has the structure $\tilde{F}_1=[\ba{cccc|c} f_{1,1} & f_{1,2} & \ldots & f_{1,{\eta_1}} & \star \ea]$, so that the control that assigns a certain $v_{1,i}$ as closed-loop eigenvector is $f_{1,i}\,x_1(t)$. 
The same procedure can be applied for $V_2$ so as to obtain a feedback matrix $\tilde{F}_2=[\ba{cccc|c} f_{2,1} & f_{2,2} & \ldots & f_{2,{\eta_2}} & \star \ea]$ in a basis adapted to $v_{2,1},\ldots,v_{2,\eta_2}$  of closed-loop eigenvectors of $A+B\,F_2\,|\,\gV_2$.
We extract $s=\dim (\gV_2+\gV_1)-\dim \gV_1$ vectors $v_{2,\alpha_1},\ldots,v_{2,\alpha_s}$ from $v_{2,1},\ldots,v_{2,\eta_2}$ (so that $\alpha_1,\ldots,\alpha_s\in \{1,\ldots,\eta_2\}$) to form a linearly independent set $v_{1,1},\ldots,v_{1,\eta_1},v_{2,\alpha_1},\ldots,v_{2,\alpha_s}$ which is a basis for an output-nulling subspace $\gV_3$ containing $\gV_1$ and $\gV_2$. In fact, we now construct a corresponding friend $\tilde{F}_3=[\ba{cccc|cccc|c} f_{1,1} & f_{1,2} & \ldots & f_{1,{\eta_1}} & f_{2,\alpha_1} & f_{2,\alpha_2} & \ldots & f_{2,{\alpha_s}} & \star \ea]$. The closed-loop system in this basis becomes
\beann
\bmat{ccc} A^F_{1,1} & 0 &  \star \\[-1mm] 0 & A^F_{s} & \star \\[-1mm] 0 & 0 &  \star \emat\,\bmat{cc} I & 0 \\[-1mm] 0 & I  \\[-1mm] 0 & 0 \emat=\bmat{cc} I & 0 \\[-1mm] 0 & I  \\[-1mm] 0 & 0 \emat\,\operatorname{diag}\{\underbrace{\l_1,\ldots,\l_1}_{\nu_1},\ldots,
\underbrace{\l_h,\ldots,\l_h}_{\nu_h},\underbrace{\l_1,\ldots,\l_1}_{\mu_1},\ldots,\underbrace{\l_h,\ldots,\l_h}_{\mu_h}\},
\eeann
where $\nu_1+\ldots+\nu_h+\mu_1+\ldots+\mu_h=\dim (\gV_1+\gV_2)$ and $\mu_1+\ldots+\mu_h=s$. It follows that $\gV_3\in \mathfrak{T}$.

We prove the second point. Let $\gK_h=\max \mathfrak{T}$. 
Clearly, $\gK_h\supseteq \ima [\ba{cccc}  V_1 & V_2 & \ldots & V_h \ea]$. The inclusion $\gK_h\subseteq \ima [\ba{cccc}  V_1 & V_2 & \ldots & V_h \ea]$ follows directly from the exhaustiveness of the parameterization of the friends of output nulling subspaces, see Section~\ref{sec:eap}.
%
%
\endproof

Theorem~\ref{lattice} showed that the largest output nulling subspace corresponding to the assignment of a certain closed-loop spectrum is -- at least in the diagonalizable case -- invariant with respect to the multiplicities. 

We now study the case when $h\ge \rho_2$; 
Lemma~\ref{lem:lemma2} can be generalized as follows. Let
\beann
&& \hspace{-0.1cm} \gT_2=\big\{q\in \NN\,|\,\exists\,F\in \mathfrak{F}(\gR^\star):\; \text{$A+B\,F\,|\,\gR^\star$ is diagonalizable and $\operatorname{card}\bigl(\sigma (A+B\,F\,|\,\gR^\star)\bigr)=q$}\big\}.
\eeann

\begin{lemma}
	\label{lem:lemma2bis}
	Let $h=\min\{\ell\in \mathbb{N}\,|\,\gR^\star=\gV^\star\cap \gS_{\ell}\}$. Then, $h=\min \gT_2$.
	\end{lemma}

\proof Using the fact that the parameterization of friends of $\gR^\star$ is exhaustive (see \cite[Theorem 3.1]{NS-SICON-14}), the proof follows along the same lines of the proof of Lemma~\ref{lem:lemma2}.\endproof

%

Loosely, for $\l$ that varies in $\real\setminus {\mathcal Z}$ the column-space of the first $n$ coordinates of a basis matrix for the kernel of $P_{\l}$ can be intuitively viewed as a subspace that ``rotates'' in $\gX$. Therefore, once the number $h$ of the closed-loop eigenvalues $\l_1,\ldots,\l_h$ is assigned, the column-space of $[\ba{cccc}  V_1 & V_2 & \ldots & V_h \ea]$ can be characterized in terms of its dimension and in terms of a subspace of it which exhibits some $\l$-invariance properties. 
The characterization in terms of the dimension was given in Theorem~\ref{th2}. 
In the next section, 
we show that the aforementioned invariant is indeed the reachability subspace on each output-nulling subspace $\ima [\ba{cccc}  V_1 & V_2 & \ldots & V_h \ea]$ obtained with different sets of closed-loop eigenvalues. 

\section{Reachability}
We now focus on the characterization of the reachability subspace on the output-nulling subspaces obtained by joining bases of the kernels of the system matrix (and, as we will recover as a particular case where $C$ and $D$ are empty, of the reachability matrix pencil).


We first introduce the following lemma. 

\begin{lemma}
	\label{lemma:diag}
	Consider a pair $(\Delta,H)$, where $\Delta $ is a $n\times n$ diagonal matrix. Let $\sigma(\Delta)=\{\l_1\,\ldots,\l_h\}$ and  $\Delta=\diag(\l_1 I_{\l_1},\ldots,\l_h I_{\l_h})$, where the matrices $I_{\l_i}$ are identity matrices of appropriate orders. Then 
	\[
	\langle \Delta\,|\im H\rangle=\im [\begin{array}{cccc}  H & \Delta\,H&\cdots &\Delta^{h-1}\,H \end{array}].
	\]
\end{lemma}
\proof
We partition $H$ conformably with $\Delta$, i.e.
$H=[\ba{ccccc} 
H_{1}^\top &  
H_{2}^\top & \ldots &
H_{h}^\top\ea]^\top$
and we show that there exist $k_1,\ldots,k_h \in \real$ such that 
\[
\Delta^h\,H=k_1\,H+k_2\,\Delta\,H+\ldots+k_{h}\,\Delta^{h-1}.
\]
The previous equation is easily seen to be equivalent to
\beann
\bmat{c}   \l_1^h\,H_{1} \\[-1.5mm]  \l_2^h\,H_{2} \\[-1.5mm] \vdots\\[-1.5mm] \l_h^h\,H_{h} \emat=
k_1\,\bmat{c} H_{1}  \\[-1.5mm] H_{2} \\[-1.5mm] \vdots \\[-1.5mm] H_{h} \emat+ 
k_2\,\bmat{c} \l_1\,H_{1}  \\[-1.5mm]  \l_2\,H_{2}   \\[-1.5mm]   \vdots   \\[-1.5mm]   \l_h\,H_{h}   \emat+ \ldots+
k_h\,\bmat{c}   \l_1^{h-1}\,H_{1}   \\[-1.5mm]   \l_2^{h-1}\,H_{2}   \\[-1.5mm]   \vdots   \\[-1.5mm]   \l_h^{h-1}\,H_{h}   \emat, 
\eeann
which can be solved in $k_i$ for every $H_{i}$ since the equation
\beann
\bmat{c} \l_1^h \\[-1.0mm]  \l_2^h\\[-1.0mm] \vdots\\[-1.0mm] \l_h^h\emat=
\bmat{ccccc} 
	1 & \l_1 & \l_1^2 & \ldots & \l_1^{h-1} \\[-1.0mm]
	1 & \l_2 & \l_2^2 & \ldots & \l_2^{h-1} \\[-1.0mm]
	\vdots & \vdots & \vdots & \ddots & \vdots \\[-1.0mm]
	1 & \l_h & \l_h^2 & \ldots & \l_h^{h-1} \emat
\bmat{c} k_1 \\[-1.0mm] k_2 \\[-1.0mm] \vdots \\[-1.0mm] k_h \emat
\eeann
is always solvable in $k_i$: in fact, the $h \times h$ matrix in the right hand-side of the latter is a Vandermonde matrix with distinct values $\l_1,\ldots,\l_h$, and it is therefore invertible.
\endproof

\begin{lemma}
	\label{lemma:reach}
	The subspace $\gV^\star_{\scriptscriptstyle \gS_h}$ (the supremal output-nulling subspace contained in $\gS_h$) is a reachability subspace.
\end{lemma}
\proof	For the sake of argument, we consider the discrete-time case. The adaptation to the continuous-time case can be obtained using e.g. the argument based on distributions of \cite[Thm. 8.22]{Trentelman-SH-01}.
This is shown in three steps. First, we show that given an arbitrary point $x_h\in\gV^\star_{\scriptscriptstyle \gS_h}$, we can always reach that point from the origin by maintaining the output at zero. Second, we show that we can always force the trajectory from $x_h$ to evolve on $\gV^\star_{\scriptscriptstyle \gS_h}$ maintaining the output at zero. Third, the trajectory between the origin and $x_h$ is entirely contained in $\gV^\star_{\scriptscriptstyle \gS_h}$. 

Since $\gV^\star_{\scriptscriptstyle \gS_h}\subseteq\gS_h$,  we have $x_h\in\gS_h$, and therefore there exist controls $u_0,\ldots,u_{h-1}$ that bring the state from the origin to $x_h$ by maintaining the output at zero: 
	\beann
		x_1 \ns&\ns = \ns&\ns B\,u_0,\quad  D\,u_0=0,\\
		x_2\ns&\ns = \ns&\ns(A+B\,F)\,x_1+B\,u_1, \quad  (C+D\,F)\,x_1+D\,u_1=0,\\
		&\vdots&\\
		x_h\ns&\ns = \ns&\ns(A+B\,F)\,x_{h-1}+B\,u_{h-1}, \quad  (C+D\,F)\,x_{h-1}+D\,u_{h-1}=0,
		\eeann
and satisfying the inclusions $x_i\in\gS_i\subseteq\gS_h$ for $i=\{1,\ldots,h-1\}$, so that the entire trajectory is in $\gS_h$.
The second point follows directly from the fact that 
at step $h$, the vector $x_h$ lies on the output-nulling subspace $\gV^\star_{\scriptscriptstyle  \gS_h}$, which implies that we can  find a control that, with initial state $x_h$, maintains the future state trajectory on $\gV^\star_{\scriptscriptstyle \gS_h}$ and the output at zero. 
We prove the third point by contradiction. Suppose that the trajectory from $0$ to $x_h$ leaves $\gV^\star_{\scriptscriptstyle \gS_h}$ (remaining in $\gS_h$ as noted above). This implies that there exists an output-nulling subspace contained in $\gS_h$ larger than $\gV^\star_{\scriptscriptstyle \gS_h}$.
\endproof

We are now ready to prove the main result of this section.


\begin{theorem}
\label{thlast}
Let $\l_1,\ldots,\l_h$ be self-conjugate, distinct and disjoint from the invariant zeros. Let
$\bsmat V_i \\[1mm] W_i\esmat$ be a basis matrix of $\ker P_{\l_i}$. Then, 
the reachability subspace
$\gR_h$ on $\ima [\ba{cccc} 
    V_1   &   V_2   &   \ldots   &   V_h   \ea]$ is 
 $\gV^\star_{\scriptscriptstyle \gS_h}$. 
\end{theorem}
\proof
 First, notice that when $h\ge \rho_2$ the statement reduces to that of Lemma~\ref{lem:lemma2bis}, and we recover the well-known result of \cite[Prop.~3]{Moore-L-78}. We now consider the case where $h<\rho_2$. Let
 \[
 \gK_h = \ima [\ba{cccc}
   V_1  &  V_2  &  \ldots  &  V_h  \ea].
 \]
 We prove that $\gR_h\subseteq\gV^\star_{\scriptscriptstyle \gS_h}$. Consider the change of coordinate matrix $T=[\ba{ccccc}    T_1  &  T_2  &  \cdots  &  T_h  &  \tilde{T}  \ea]$ in $\gX$ such that
$\im [\ba{ccccc}    T_1  & \ldots  &  T_k  \ea]=\im V_1+\cdots+\im V_k$ for all $k\in \{1,\ldots,h\}$.
Consider an $m\times m$ non-singular matrix $\Omega=[\ba{cc}   \Omega_1  &  \Omega_2  \ea]$ such that $\Omega_1$ is a basis matrix
for $B^{-1}\gK_h\cap \ker D$. 
By construction, 
\[
T^{-1}(A+B\,F)T=\bmat{cc} A^F_{1,1} & \star\\[-1.0mm]
0 & \star\emat,\qquad T^{-1}B\,\Omega=\bmat{cc} B^{\Omega}_{1,1} & \star \\[-1.0mm] 0  & \star \emat.
\]
Since the closed-loop map restricted to $\gK_h$, i.e. $A^F_{1,1}$, is diagonalizable, there exists a  nonsingular matrix $\hat{S}=\bsmat S && 0\\[1mm] 0 && I \esmat$ such that $S^{-1}A^F_{1,1}S=\Delta$, where $\Delta$ is diagonal, and 
\[
\hat{S}^{-1}T^{-1}(A+B\,F)T\,\hat{S}=\bmat{cc} \Delta & \star\\[-1.0mm]
0 & \star\emat,\quad \hat{S}^{-1}T^{-1}B\,\Omega=\bmat{cc} H & \star \\[-1.0mm] 0  & \star \emat
\]
for a suitable submatrix $H$.
In the new coordinates, the reachable subspace on $\gK_h$ can be computed as $\gR_h = \im \bsmat H && \Delta\,H && \cdots && \Delta^{\e-1}\,H\\[1mm] 0 && 0 && 0 && 0\esmat$, for a sufficiently large $\e$. In view of Lemma \ref{lemma:diag}, we have $\e=h$, which implies that for every point $x_h\in \gR_h$ there exist $u_0,\ldots,u_{h-1}$ such that, by recalling that $\gR_h$ is output-nulling, there hold 
\beann
x_1\ns&\ns=\ns&\ns B\,u_0,\quad  D\,u_0=0,\\
x_2\ns&\ns=\ns&\ns (A+B\,F)\,x_1+B\,u_1, \quad  (C+D\,F)\,x_1+D\,u_1=0,\\
&\vdots&\\
x_h\ns&\ns=\ns&\ns (A+B\,F)\,x_{h-1}+B\,u_{h-1}, \quad  (C+D\,F)\,x_{h-1}+D\,u_{h-1}=0.
\eeann
The previous equalities show that $x_h\in\gS_h$, which, considering that $\gV^\star_{\scriptscriptstyle \gS_h}$ is the largest output nulling on $\gS_h$, imply that $x_h \in \gV^\star_{\scriptscriptstyle \gS_h}$, from which we have $\gR_h\subseteq\gV^\star_{\scriptscriptstyle \gS_h}$.

Now we prove that $\gV^\star_{\scriptscriptstyle \gS_h}\subseteq\gK_h$. Given a change of coordinates matrix
$T=[\ba{cc}  T_1 &   T_2  \ea]$, with $\im T_1=\gV^\star_{\scriptscriptstyle \gS_h}$ and an $m\times m$ non-singular matrix $\Omega=[\ba{cc}   \Omega_1  &  \Omega_2  \ea]$ such that $\Omega_1$ is a basis matrix
for $B^{-1}\gV^\star_{\scriptscriptstyle \gS_h}\cap \ker D$, for every friend $F$ of $\gV^\star_{\scriptscriptstyle \gS_h}$, we obtain
\[
T^{-1}(A+B\,F)T=\bmat{cc} A_{11} & \star\\[-1.0mm]
0 & \star\emat,\qquad T^{-1}B\,\Omega=\bmat{cc} B_{11} & \star \\[-1.0mm] 0 & \star \emat
\]
In view of Lemma~\ref{lemma:reach}, $\gV^\star_{\scriptscriptstyle \gS_h}$ is a reachability output-nulling subspace, and therefore the pair $(A_{11},B_{11})$ is reachable and, in the new basis, $\gV^\star_{\scriptscriptstyle \gS_h}=\ima \bsmat I \\[1mm] 0 \esmat = \ima \bsmat B_{11} && A_{11}B_{11} && \cdots && A_{11}^{h-1} B_{11} \\[1mm] 0 && 0 && \ldots && 0 \esmat$
because $\gV^\star_{\scriptscriptstyle \gS_h}\subseteq\gS_h$, so that each point can be reached in at most $h$ iterations, in the discrete case, or using the same argument based on distributions (Dirac deltas and their distributional derivatives) in the continuous time, see \cite[Thm. 8.22]{Trentelman-SH-01}. Moreover, $\gV^\star_{\scriptscriptstyle \gS_h} \supset \ima \bsmat B_{11} && A_{11}B_{11} && \cdots && A_{11}^{h-2} B_{11} \\[1mm] 0 && 0 && \ldots && 0 \esmat$ because $\gS_{h-1} \subset \gS_h$ and $\gsV \cap \gS_j $ is nondecreasing in $j$ and converges to $\gsR$. Therefore, in view of Lemma~\ref{lem:lemma2} adapted to $(A_{11},B_{11})$ (instead of $(A,B)$), $h$ is the minimum number such that there exists $F_{11}$ that can assign $\sigma(A_{11}+B_{11}\,F_{11})=\{\l_1,\ldots,\l_h\}$ 
in such a way that $A_{11}+B_{11}\,F_{11}$ is diagonalizable.
 Since the set of all friends of $\gV^\star_{\scriptscriptstyle \gS_h}$ can be parameterized as
$F + \Omega \bsmat F_{11} && \star \\[1mm] 0 && \star \esmat T^{-1}$ where $F_{11}$ as well as the entries indicated with $\star$ are arbitrary (see \cite[Thm 2]{NP-17}), there exists a friend $\widehat F$ of $\gV^\star_{\scriptscriptstyle \gS_h}$ that can assign $\sigma (A+B\,\widehat F\,|\,\gV^\star_{\scriptscriptstyle \gS_h})=\{\l_1,\ldots,\l_h\}$ 
in such a way that $A+B\,\widehat F\,|\,\gV^\star_{\scriptscriptstyle \gS_h}$ is diagonalizable.
Recalling that $\gK_h$ is the largest output nulling subspace where $\{\l_1,\ldots,\l_h\}$ can be assigned without Jordan forms, we have $\gV^\star_{\scriptscriptstyle \gS_h}\subseteq\gK_h$. Finally, the fact that $\gV^\star_{\scriptscriptstyle \gS_h}$ is reachable together with the fact that $\gR_h$ is the reachability subspace on $\gK_h$ implies that $\gV^\star_{\scriptscriptstyle \gS_h}\subseteq\gR_h$ and, consequently that $\gV^\star_{\scriptscriptstyle \gS_h}=\gR_h$.
\endproof

Notice that as a result of this theorem, $\gK_h\cap B\ker D=\gV^\star_{\scriptscriptstyle \gS_h}\cap B\ker D$ for all possible sets $\{\l_1,\ldots,\l_h\}$.

\begin{remark}{
When $h\ge \rho_2$, the statements of Theorem~\ref{th2}
and Theorem~\ref{thlast} allow us to 
to write the chain of identities
\[
\gR^\star=\gK_h=\gV^\star \cap \gS_h=\gR_h=\gV^\star_{\gS_h}=\gV^\star_{\gS^\star}=\gV^\star\cap \gS^\star,
\]
and we recover the well-known identity $\gR^\star=\gV^\star\cap \gS^\star$ \cite{Morse-73-SIAM}. We also note that when $h=\rho_2$, we have $\gS_h=\gS^\star$, so that the sequence to generate the infimal input-containing $\gS^\star$ becomes stationary at the step at which adding a new kernel the image does not change.
}
\end{remark}

%
%

We now consider the case where the kernels are extracted from the reachability matrix pencil instead of the system matrix pencil.


\begin{corollary}
\label{last}
Let $\l_1,\ldots,\l_h$ be disjoint from the uncontrollable eigenvalues. Let
$\bsmat V_i \\[1mm] W_i\esmat$ be a basis matrix of $\ker S_{\l_i}$.
Then, the reachable subspace $\gR_h$ on $\ima [\ba{cccc}
 V_1 & V_2 & \ldots & V_h\ea]$ is the largest $(A,B)$-controlled invariant contained in $\ima [\ba{cccc}  B & A\,B & \ldots & A^{h-1}B\ea]$.
\end{corollary}
\proof
When $C$ and $D$ are empty, $\gS_h$ is the reachable subspace in $h$ steps. In other words, $\gS_h=\ima [\ba{cccc} B & A\,B & \ldots & A^{h-1}B\ea]$, and the uncontrollable eigenvalues are the invariant zeros of the system. The rest of the proof follows as a particular case of the proof of Theorem~\ref{thlast}.
\endproof

The reachable subspace $\gR_h$ of Corollary~\ref{last} can be computed using the subspace sequence that converges to $\gV^\star_{\gE}$ by considering  $C$ and $D$ to be empty and $\gE=\ima [\ba{cccc} B & A\,B & \ldots & A^{h-1}B\ea]$.

\end{document}